\documentclass[11pt]{article}

\usepackage{amsmath,amssymb}

\usepackage{times}
\usepackage{bm}
\usepackage{algorithm}
\usepackage{xargs}[2008/03/08]
\usepackage{amssymb}
\usepackage{mathrsfs}
\usepackage{graphicx}
\usepackage{rotating,subfigure}
\usepackage[flushleft]{threeparttable}
\usepackage{multirow}
\usepackage{enumitem} 

\usepackage{star}

\usepackage[colorlinks,
linkcolor=blue,
anchorcolor=blue,
citecolor=blue
]{hyperref}

\def\##1\#{\begin{align}#1\end{align}}
\def\$#1\${\begin{align*}#1\end{align*}}

\newcommand{\Rom}[1]{\text{\uppercase\expandafter{\romannumeral #1\relax}}}

\usepackage{txfonts}

\usepackage{geometry}
\numberwithin{equation}{section}

\allowdisplaybreaks[4]

\begin{document}

\title{ \LARGE Limiting laws for extreme eigenvalues of large-dimensional spiked Fisher matrices with a divergent number of spikes}    

\author{Junshan Xie\thanks{Co-first author. School of Mathematics and Statistics, Henan University, Kaifeng, China.},~  Yicheng Zeng\thanks{Co-first author. Department of Mathematics, Hong Kong Baptist University, Hong Kong.},~  Lixing Zhu\thanks{Corresponding author. Research Center for Statistics and Data Science, Beijing Normal University, Zhuhai, China and Department of Mathematics, Hong Kong Baptist University, Hong Kong. Email address: \url{lzhu@hkbu.edu.hk}.}}

\date{ }

\maketitle

\vspace{-0.25in}

\begin{abstract}
Consider the $p\times p$  matrix that is the product of a population covariance matrix and the inverse of another population covariance matrix. Suppose that their difference has a divergent rank with respect to $p$,  when two samples of sizes $n$ and $T$ from the two populations are available, we construct its corresponding sample version. In the regime of high dimension where both $n$ and $T$ are proportional to $p$, we investigate the limiting laws for extreme (spiked) eigenvalues of the sample (spiked) Fisher matrix when the number of spikes is divergent and these spikes are unbounded. 
\end{abstract}

\noindent
{\bf Keywords:} Extreme eigenvalue, Fisher matrix, Phase transition phenomenon, Random matrix theory, Spiked population model.

\section{Introduction}\label{sec1}
In the last few decades, as the remarkable development in storage devices and computing capability, the demand for processing complex-structured data increases dramatically.  One of the features as well as the challenges of these data sets is their high dimensions. The difficulty is that the classical limit theory for multivariate statistical analysis fails to ensure reliable inference for high-dimensional data analysis. Classical limit theorems require ``small $p$ large $n$'' to keep their validity, which conflicts with the situation ``large $p$ large $n$'' in high-dimensional settings in the sense that  $p/n\rightarrow c>0$ as the asymptotic properties are rather different. To attack the relevant issues,  random matrix theory (RMT) serves as a powerful tool in addressing statistical problems in high dimensions.
The first research of random matrices in multivariate statistics was about the Wishart matrices  in \cite{wishart1928generalised}.
Abundant research has been established for various topics in this field during the past half century, especially in recent years.
In the area of RMT in statistics, we refer to monographs \cite{bai2010spectral} and \cite{yao2015large} for systematical study and \cite{debashis2014random} for a comprehensive review.

A relevant topic in multivariate statistics is about testing the equality of two covariance matrices:
\begin{align}\label{hypotheis}
H_0:\ \Sigma_1=\Sigma_2\quad \text{vs.}\quad
H_1:\ \Sigma_1=\Sigma_2+\Delta,
\end{align}
where $\Sigma_1$ and $\Sigma_2$ are two covariance matrices corresponding to two $p$-variate populations, and $\Delta$ is a non-negative definite matrix with rank $q$. Let $\mathbf{S}_1$ and $\mathbf{S}_2$ be the sample covariance matrices from these two populations, respectively.
When $\textbf{S}_2$ is invertible, the random matrix $\textbf{F}= \textbf{S}_2^{-1}\textbf{S}_1$ is called a Fisher matrix.

The difference between the null hypothesis and the alternative hypothesis relies on those extreme eigenvalues of $\mathbf{F}$.
Under the null hypothesis,  $\Sigma_1=\Sigma_2$, \cite{wachter1980limiting} established the well-known Wacheter distribution as the limiting spectral distribution (LSD) of $\mathbf{F}$. Some extensions were built later (see examples in \cite{silverstein1985limiting}, \cite{silverstein1995strong} and \cite{silverstein1995empirical}). Furthermore, \cite{bai1998no} pointed out the fact that the largest eigenvalue of $\mathbf{F}$ converges to the upper bound of the support of the LSD of $\mathbf{F}$.
Under the alternative hypothesis, $\mathbf{F}$ is called a {\it spiked Fisher matrix} (see \cite{wang2017extreme}), because $\Sigma_2^{-1}\Sigma_1$ has a spiked structure similar to that of a {\it spiked population model} proposed by \cite{johnstone2001distribution}. More specifically, the matrix $\Sigma_2^{-1}\Sigma_1$ is assumed to have the spectrum
\begin{align}\label{spectrum}
\text{spec}(\Sigma_2^{-1}\Sigma_1\mathbf)=\{\lambda_1,\ldots, \lambda_q,1,\ldots,1\},
\end{align}
where $\lambda_1\ge\ldots\ge\lambda_q>1$.
When the rank $q$ of $\Delta$ is finite, \cite{dharmawansa2014local} showed the phase transition phenomenon of the extreme eigenvalues of $\mathbf{F}$ under Gaussian population assumption. That is, for $1\le i\le q$, the $i$-th largest eigenvalue of $\mathbf{F}$ will depart from the upper bound of the support of LSD of $\mathbf{F}$ if and only if $\lambda_1$ exceeds certein phase transition point.  \cite{wang2017extreme} extended it to the cases without Gaussian assumption and established central limit theorems for the outlier eigenvalues of $\mathbf{F}$.

We in this paper consider, as a reasonable extension in theory and applications, the case of divergent $q$ with respect to the dimension $p$. We will investigate the convergence in probability and central limit theorems for spiked eigenvalues of spiked Fisher matrices. We formulate our problem as follows.

Assume that
\begin{align}\label{matrix YZ}
\mathbf{Y}=\left(\mathbf{y}_{1},\ldots,\mathbf{y}_{T}\right)=(y_{ij})_{1\leq i\leq p,1\leq j\leq T}\in \mathbb{R}^{p\times T}\quad \text{and}\quad
\mathbf{Z}=\left(\mathbf{z}_{1},\ldots,\mathbf{z}_{n}\right)=(z_{ij})_{1\leq i\leq p,1\leq j\leq n}\in \mathbb{R}^{p\times n}
\end{align}
are two independent arrays of independent real-valued random variables with zero mean and unit variance. We consider two samples $\{\Sigma_1^{1/2}\mathbf{y}_{i}\}_{1\le i\le T}$ and $\{\Sigma_2^{1/2}\mathbf{z}_{i}\}_{1\le i\le n}$, then their corresponding sample covariance matrices can respectively be written as
\begin{align*}
\textbf{S}_{1}=\frac{1}{T}\sum_{i=1}^{T}\Sigma_1^{\frac{1}{2}}\mathbf{y}_{i}\mathbf{y}_{i}^{\top}\Sigma_1^{\frac{1}{2}}=\frac{1}{T}\Sigma_1^{\frac{1}{2}}\mathbf{Y}\mathbf{Y}^{\top}\Sigma_1^{\frac{1}{2}}\quad \text{and}\quad
\textbf{S}_{2}=\frac{1}{n}\sum_{i=1}^{n}\Sigma_2^{\frac{1}{2}}\mathbf{z}_{i}\mathbf{z}_{i}^{\top}\Sigma_2^{\frac{1}{2}}=\frac{1}{n}\Sigma_2^{\frac{1}{2}}\mathbf{Z}\mathbf{Z}^{\top}\Sigma_2^{\frac{1}{2}}.
\end{align*}
Also,  define the Fisher matrix
$
 \textbf{F}:=\textbf{S}_{2}^{-1}\textbf{S}_{1},
$
as the sample version of matrix $\Sigma^{-1}_2\Sigma_1$.
We aim to investigate the limiting properties of the eigenvalues of $\textbf{F}$. As the eigenvalues of $\textbf{F}$ remain invariant under the linear transformation
\begin{align}\label{linear_transformation}
\left(\textbf{S}_{1},\textbf{S}_{2}\right)\rightarrow \left(\Sigma_{2}^{-\frac{1}{2}} \textbf{S}_{1}\Sigma_{2}^{-\frac{1}{2}} ,\Sigma_{2}^{-\frac{1}{2}} \textbf{S}_{2}\Sigma_{2}^{-\frac{1}{2}}\right),
\end{align}
thus   we can assume $\Sigma_{2}=\textbf{I}_{p}$ throughout this paper without loss of generality. Under the assumption~\eqref{spectrum}, eigenvalues of $\Sigma_1$ are $\lambda_1\ge \ldots\ge \lambda_q>\lambda_{q+1}=\ldots=\lambda_p=1$.
Recalling \eqref{hypotheis} that $\Sigma_1$ is a rank $q$ pertubation of $\Sigma_2=\mathbf{I}_p$, we simply assume
\begin{align}\label{Sigma1}
\Sigma_1=\begin{pmatrix}
\Sigma_{11}&0\\
0&\mathbf{I}_{p-q}
\end{pmatrix}.
\end{align}
For the sake of brevity and readability,  we write the eigenvalues of $\textbf{F}$ in descending order  $\hat\lambda_{1}\geq  \ldots\geq \hat\lambda_{p}$, simplifying the double subscripts as single ones. It should be noted that $\hat\lambda_i$ is related to the sample size $n$.

We then describe the related work and our contributions in this paper.
When the number of the spiked eigenvalues $q$ is fixed, and all the spiked eigenvalue $\lambda_{i},i=1,\ldots,q$, are bounded, there are some results on the limiting properties of the eigenvalues of $\textbf{F}$ in the literature. Such as, the almost surely convergence (strong consistency)  and central limit theorem (CLT) of spiked eigenvalues (\cite{wang2017extreme}) and asymptotically Tracy-Widom distribution for the largest non-spiked eigenvalue (\cite{han2018unified} and \cite{han2016tracy}).
In this paper, we consider the case that the number of spiked eigenvalues $q=q(p)\rightarrow\infty$ as $p\rightarrow\infty$, and spiked eigenvalues $\lambda_{i},1\leq i \leq q$ diverge as $p\rightarrow\infty$.
To the  best of our knowledge, there is no relevant result in the literature. A relevant work is \cite{cai2017limiting} who studied spiked population models, where the asymptotics for spiked eigenvalues, including convergence in probability (weak consistency) and CLT, as well as Tracy-Widom law for the largest nonspiked eigenvalue were built under a quite general framework.
Unlike the case of fixed $q$ and bounded spikes $\lambda_i$, $1\le i\le q$, normalizations for $\hat\lambda_i$, $1\le i\le q$ are needed for the divergent $q$ case.
Consider the normalized eigenvalues $\hat\lambda_i/\lambda_i$ in consistency and $(\hat\lambda_i-\theta_i)/\theta_i$ in CLT, where $\theta_i$ is a centered parameter defined later.

The basic approach behind the proofs of the asymptotics for spiked eigenvalues is the analysis of an equation for the determinant of a $q\times q$ random matrix (indexed by $n$). When $q$ is bounded, \cite{wang2017extreme} derived the almost sure entrywise convergence of the $q\times q$ matrix (and hence the convergence with respect to matrix norms) and then solving the equation to lead to the almost sure limits of spiked eigenvalues. This argument does not work in the divergent $q$ case where the convergence of a $q\times q$ matrix with respect to some norm could not be directly implied by the entrywise convergence. Instead, we use the CLT for random sequilinear forms in \cite{bai2008central} to derive the convergence rate of each entry, and then use Chebyshev's inequality to put all entries together to derive the convergence rate of the matrix in $\ell_\infty$ norm. In this way, we achieve the convergence in probability as well as the CLT of spiked eigenvalues (after proper normalizations). This approach is similar to that used in \cite{cai2017limiting}, so some technical assumptions are also imposed similarly.

The remaining parts of the paper are organized as follows. Section \ref{sec2} establishes the main results, including the convergence in probability of $\hat{\lambda}_i/\lambda_i$ and central limit theorems of $(\hat\lambda_i-\theta_i)/\theta_i$, for those spiked eigenvalues of spiked Fisher matrix $\mathbf{F}$. Here, $\theta_i$, $1\le i\le q$, is a sequence of centering parameters defined in this section. In the Section~\ref{sec3}, we show the proofs of our main results in Section~\ref{sec2}. Some important technical lemmas and their proofs are displayed in the Section~\ref{sec4}.

\section{Main results}\label{sec2}
\subsection{Notations and assumptions}
Considering the linear transformation \eqref{linear_transformation}, we assume that $\Sigma_{2}=\textbf{I}_{p}$ without loss of generality, and then $\Sigma_1$ has the structure as shown in \eqref{Sigma1}. Further, we decompose the $\Sigma_{11}$ in \eqref{Sigma1} as
\begin{align*}
\Sigma_{11}=\mathbf{U}^{\top}\Lambda_{1} \mathbf{U}.
\end{align*}
Here, $\mathbf{U}\equiv(\mathbf{u}_1,\mathbf{u}_2,\ldots,\mathbf{u}_q)^{\top}$ is a $q\times q$ orthogonal matrix and
\begin{align*}
\Lambda_{1}=\text{diag}(\underbrace{\lambda_{1},\ldots,\lambda_{N_1}}_{n_{1}},\ldots,\underbrace{\lambda_{N_{\ell-1}+1},\ldots,\lambda_{q}}_{n_{\ell}}),
\end{align*}
where $\lambda_1=\ldots=\lambda_{N_1}>\ldots>\lambda_{N_{\ell-1}+1}=\ldots=\lambda_{q}$ and $N_i:=\sum_{j=1}^i n_j$ for $1\le i\le \ell$. In this case, $\Sigma_1$ can be decomposed as
\begin{align*}
\Sigma_1=\begin{pmatrix}
\mathbf{U}^\top&0\\
0&\mathbf{I}_{p-q}
\end{pmatrix}
\begin{pmatrix}
\Lambda_1&0\\
0&\mathbf{I}_{p-q}
\end{pmatrix}
\begin{pmatrix}
\mathbf{U}&0\\
0&\mathbf{I}_{p-q}
\end{pmatrix}
=:
\begin{pmatrix}
\mathbf{U}^\top&0\\
0&\mathbf{I}_{p-q}
\end{pmatrix}
\Lambda
\begin{pmatrix}
\mathbf{U}&0\\
0&\mathbf{I}_{p-q}
\end{pmatrix}.
\end{align*}

We give  decompositions of the sample covariance matrices $\textbf{S}_{1}$ and $\textbf{S}_{2}$ as follows. We fisrt decompose the matrices $\mathbf{Y}$  and $\mathbf{Z}$ defined in \eqref{matrix YZ} as $\mathbf{Y}=( \mathbf{Y}_{1}^{\top}, \mathbf{Y}_{2}^{\top}
)^{\top}$ and $\mathbf{Z}=( \mathbf{Z}_{1}^{\top}, \mathbf{Z}_{2}^{\top}
)^{\top}$, where $\mathbf{Y}_{1}, \mathbf{Z}_{1}\in \mathbb{R}^{q\times n}$ and $\mathbf{Y}_{2},\mathbf{Z}_{2}\in \mathbb{R}^{(p-q)\times n}$.
Let $\mathbf{X}:=\Sigma_1^{1/2} \mathbf{Y}$. Then we can similarly write $\mathbf{X}=( \mathbf{X}_{1}^{\top}, \mathbf{X}_{2}^{\top}
)^{\top}$, where $\mathbf{X}_{1}=\Sigma_{11}^{1/2}\mathbf{Y}_{1}=\mathbf{U}^{\top}\Lambda_{1}^{1/2}\mathbf{U}\mathbf{Y}_{1}\in \mathbb{R}^{q\times T}$ and $\mathbf{X}_{2}=\mathbf{Y}_2\in \mathbb{R}^{(p-q)\times T}$. It follows that
\begin{align}\label{decomposition}
\textbf{S}_{1}=\begin{pmatrix} \frac{1}{T}\mathbf{X}_{1}\mathbf{X}_{1}^{\top}&\frac{1}{T}\mathbf{X}_{1}\mathbf{X}_{2}^{\top}\\\frac{1}{T}\mathbf{X}_{2}\mathbf{X}_{1}^{\top}&\frac{1}{T}\mathbf{X}_{2}\mathbf{X}_{2}^{\top}\end{pmatrix}
\quad \text{and}\quad
\textbf{S}_{2}=\begin{pmatrix} \frac{1}{n}\mathbf{Z}_{1}\mathbf{Z}_{1}^{\top}&\frac{1}{n}\mathbf{Z}_{1}\mathbf{Z}_{2}^{\top}\\\frac{1}{n}\mathbf{Z}_{2}\mathbf{Z}_{1}^{\top}&\frac{1}{n}\mathbf{Z}_{2}\mathbf{Z}_{2}^{\top}\end{pmatrix}.
\end{align}

For $\lambda\in \mathbb{R}\setminus \{0\}$, we introduce
\begin{gather}
\textbf{F}_{0}=\left(\frac{1}{n}\mathbf{Z}_{2}\mathbf{Z}_{2}^{\top}\right)^{-1}\left(\frac{1}{T}\mathbf{X}_{2}\mathbf{X}_{2}^{\top}\right),\nonumber\quad 
\mathbf{M}(\lambda)= \textbf{I}_{p-q}- \frac{\textbf{F}_{0}}{\lambda},\nonumber \\
\widetilde{ m}_{\theta}(z)= \frac{1}{p-q}\text{tr} \left(z\textbf{I}_{p-q}-\frac{\textbf{F}_{0}}{\theta} \right)^{-1},\ \theta\in \mathbb{R},\ z\in \mathbb{C}^{+}.\label{def_m}
\end{gather}
Let $\mu_{1}\geq \ldots \ge\mu_{p-q}$  be the eigenvalues of the Fisher matrix $\textbf{F}_{0}$.  Then the empirical spectral distribution (ESD) of $\textbf{F}_{0}$ can be defined as
\begin{align*}
F_{n}(x)=\frac{1}{p-q}\sum_{j=1}^{p-q}\textbf{1}_{\{\mu_{j}\leq x\}},\ x\in \mathbb{R}.
\end{align*}
By the result in \cite{wang2017extreme}, under the assumption of $p/n\rightarrow y\in(0,1)$ and $p/T\rightarrow c>0$, almost surely, the empirical spectral distribution $F_{n}$ weakly  converges to the limiting spectral distribution $F_{c,y}$, whose Stieltjes transform
$\mathcal{S}(z)=\int_{-\infty}^{\infty}(x-z)^{-1}dF_{c,y}(x)$ satisfies, for $ z\notin [a, b]$
\begin{align}\label{sz}
\mathcal{S}(z)=\frac{1-c}{zc}-\frac{ c[z(1-y)+1-c]+2zy-c\sqrt{[z(1-y)+1-c]^{2}-4z}   }{2zc(c+zy)},
\end{align}
where $a=( 1-\sqrt{c+y-cy})^2(1-y)^{-2}$ and $b=(1+\sqrt{c+y-cy})^2(1-y)^{-2}$.

In the following, for any complex matrix $A$, we use $s_i(A)$ to denote the $i$-th largest singular value, and $\|A\|$ to denote the largest singular value  throughout the paper.
Write $a_n={\rm O}_{a.s.}\left(b_n\right)$ if it almost surely holds that $a_n={\rm O}\left(b_n\right)$. Throughout this paper $C$ is a constant that may vary from place to place.

The following assumptions are required.
\begin{assumption}\label{A1}
 $y_{p}:=p/n\rightarrow y\in(0,1)$, $\tilde{y}_p:=(p-q)/n$; $c_{p}:=p/T\rightarrow c>0$, $\tilde{c}_p:=(p-q)/T$;
$q=q(n)\rightarrow\infty$ as $n\rightarrow\infty$ but $q={\rm o }(n^{\frac{1}{6}})$.
\end{assumption}
\begin{assumption}\label{A2}
For any $1\le i\le q$, $\lambda_i$ satisfies $q^2/\lambda_i\rightarrow 0$ and either of the two  following conditions:\\
$(a).$ $\lambda_i^{-1}\sum_{j=1}^q \lambda_j={\rm o}(q^{-\frac{1}{2}}n^{\frac{1}{4}})$ and $\lambda_i^{-2}\sum_{j=1}^q\lambda_j={\rm o} (q^{-1})$;\  $(b).$ $\lambda_i\sum_{j=1}^q \lambda_j^{-1}={\rm o} (q^{-\frac{1}{2}}n^{\frac{1}{4}})$.
\end{assumption}
\begin{assumption}\label{A3}
Random vectors in $\{\mathbf{y}_i: 1\le i\le T\}\bigcup\{\mathbf{z}_i: 1\le i\le n\}$ are independent identically distributed, ${\rm E}z_{ij}=0$, ${\rm E}|z_{ij}|^2=1$ $\forall 1\le i\le p,\ 1\le j\le n$ and $\sup_{1\leq i\leq p}{\rm E}|z_{ij}|^{4}<\infty$.
\end{assumption}
\begin{assumption}\label{A4}
There exists a constant $C>1$ such that $\lambda_{N_i}/\lambda_{N_{i+1}}\geq C$ for any $1\le i\le \ell -1$.
\end{assumption}
\begin{assumption}\label{A5}
Suppose that $\{\lambda_i\}_{1\le i\le q}$ are of bounded multiplicities, i.e., $\sup_{1\le i\le \ell}n_i<\infty$.
\end{assumption}

\subsection{Weak consistency}
The weak consistency of  $\hat\lambda_i$ is stated below. Due to the fact that $\lambda_i$ may go to infinity with $n$,  consider the limit in probability for the ratio $\hat\lambda_i/\lambda_i$, $1\le i\le q$.

\begin{theorem}\label{th1.2}
Assume that Assumptions~\ref{A1}, \ref{A2} and \ref{A3} hold. Then  for all $1\leq i \leq q$,
\begin{align*}
\frac{\hat \lambda_i}{\lambda_i}=\frac{1}{1-y}+{\rm O }\left(y_p-y\right)+\kappa q\cdot {\rm O}_p \left(\frac{1}{\sqrt n}+\lambda_i^{-1}\right),
\end{align*}
where $\kappa:=\min\{\kappa_1,\kappa_2\}$ with $\kappa_1:= q+\lambda_i^{-1}\sum_{j=1}^q\lambda_j$ and $\kappa_2:= q+\lambda_i\sum_{j=1}^q\lambda_j^{-1}$.
\end{theorem}
\begin{remark}
Note that the limit of the ratio $\hat{\lambda}_i/\lambda_i$  is $1/(1-y)>1$, for all $1\le i\le q$. This is different from the relevant limit  for spiked population model with divergent $q$, which is $1$ (see Theorem~2.1 in \cite{cai2017limiting}). Roughly speaking, when we take $y\rightarrow 0$ with $1/(1-y)\rightarrow 1$, asymptotically,  a spiked Fisher matrix behaves similarly as the  sample covariance matrix in a spiked population model.
\end{remark}
\begin{remark}
In the case of fixed $q$ and bounded spikes $\lambda_i$, $1\le i\le q$,
Theorem~3.1 in \cite{wang2017extreme} shows that almost surely the spiked eigenvalue
$\hat{\lambda}_i$  converges  to the limit
$\lambda_i(\lambda_i+c-1)(\lambda_i-\lambda_i y-1)^{-1}$. Simply taking
$\lambda_i\rightarrow \infty$, the limit of $\lambda_i(\lambda_i+c-1)(\lambda_i-\lambda_i y-1)^{-1}\lambda_i^{-1}$  equals to $1/(1-y)$. Thus, Theorem~\ref{th1.2} indicates that for  the divergent $q$ case, the result coincides with the result for the fixed $q$ case in \cite{wang2017extreme}.
\end{remark}
\begin{remark}
In Theorem~\ref{th1.2} we only consider unbounded spikes, but actually it can be readily extended to handle the case with both bounded and unbounded spikes. Consider the model
\begin{align*}
\Sigma_1=
\begin{pmatrix}
\mathbf{U}^\top&0\\
0&\mathbf{I}_{p-q}
\end{pmatrix}
\Lambda
\begin{pmatrix}
\mathbf{U}&0\\
0&\mathbf{I}_{p-q}
\end{pmatrix},
\end{align*}
where $\Lambda={\rm diag}(\lambda_1,\ldots,\lambda_q,\lambda_{q+1},\ldots,\lambda_{q+q_0},1,\ldots,1)$, $q={\rm o }(n^{1/6})$ and $q_0$ is bounded. Assume that spikes $\lambda_1\geq \ldots \geq \lambda_q$ are unbounded as in Theorem~\ref{th1.2} and $\lambda_{q+1}\geq \ldots \geq \lambda_{q+q_0}$ are bounded.
For $q+1\le i\le q+q_0$, by Theorem~A.10 in \cite{bai2010spectral}, we have
\begin{align*}
\hat \lambda_i= s_i\left(\mathbf{S}_2^{-1}\mathbf{S}_1\right)\le s_i\left(\mathbf{S}_1\right) s_1\left(\mathbf{S}_2^{-1}\right)\le s_i\left(\Sigma_1\right)s_1\left(\frac{1}{T}\mathbf{Y}\mathbf{Y}^\top\right)s_1\left(\mathbf{S}_2^{-1}\right)<\infty
\end{align*}
almost surely. So it holds that
\begin{align*}
\det\left(\frac{\hat\lambda_{i}}{n}\mathbf{Z}_{1}\mathbf{Z}_{1}^{\top}-\frac{1}{T}\mathbf{X}_{1}\mathbf{X}_{1}^{\top}\right)\neq 0.
\end{align*}
Similar to the decomposition in \eqref{product}, we have
{\small
\begin{gather}\label{extension1}
\det\left\{\left(\frac{\hat\lambda_{i}}{n}\mathbf{Z}_{2}\mathbf{Z}_{2}^{\top}-\frac{1}{T}\mathbf{X}_{2}\mathbf{X}_{2}^{\top}\right)
\displaystyle -\left(\frac{\hat\lambda_{i}}{n}\mathbf{Z}_{2}\mathbf{Z}_{1}^{\top}-\frac{1}{T}\mathbf{X}_{2}\mathbf{X}_{1}^{\top}\right)
\left(\frac{\hat\lambda_{i}}{n}\mathbf{Z}_{1}\mathbf{Z}_{1}^{\top}-\frac{1}{T}\mathbf{X}_{1}\mathbf{X}_{1}^{\top}\right)^{-1}
\left(\frac{\hat\lambda_{i}}{n}\mathbf{Z}_{1}\mathbf{Z}_{2}^{\top}-\frac{1}{T}\mathbf{X}_{1}\mathbf{X}_{2}^{\top}\right)  \right\}=0.
\end{gather}
}
In the same manner as used in the proof of Theorem~\ref{th1.2}, it can be checked that
\begin{align*}
\Big|\Big|\Big|\left(\frac{\hat\lambda_{i}}{n}\mathbf{Z}_{2}\mathbf{Z}_{1}^{\top}-\frac{1}{T}\mathbf{X}_{2}\mathbf{X}_{1}^{\top}\right)
\left(\frac{\hat\lambda_{i}}{n}\mathbf{Z}_{1}\mathbf{Z}_{1}^{\top}-\frac{1}{T}\mathbf{X}_{1}\mathbf{X}_{1}^{\top}\right)^{-1}
\left(\frac{\hat\lambda_{i}}{n}\mathbf{Z}_{1}\mathbf{Z}_{2}^{\top}-\frac{1}{T}\mathbf{X}_{1}\mathbf{X}_{2}^{\top}\right)\Big|\Big|\Big|_\infty={\rm o}_p(1).
\end{align*}
Then the solution of equation \eqref{extension1} is close to that of the equation
\begin{align}\label{extension2}
\det\left(\frac{\hat\lambda_{i}}{n}\mathbf{Z}_{2}\mathbf{Z}_{2}^{\top}-\frac{1}{T}\mathbf{X}_{2}\mathbf{X}_{2}^{\top}\right)=0.
\end{align}
Note that the solution of \eqref{extension2} is an eigenvlaue of the spiked Fisher matrix $(\mathbf{Z}_{2}\mathbf{Z}_{2}^{\top}/n)^{-1}(\mathbf{X}_{2}\mathbf{X}_{2}^{\top}/T)$ which has been well studied by \cite{wang2017extreme}. Thus, the weak consistency for all outliers $\hat\lambda_i$, $1\le i\le q+q_0$, could be achieved by combining Theorem~3.1 in \cite{wang2017extreme} and Theorem~\ref{th1.2}. Such a kind of extension could also be considered for the CLT in Theorem~\ref{th1.3}.
\end{remark}

\subsection{Central limit theorem}
As $\lambda_i$, $1\le i\le q$, goes to infinity, the consistency of $\hat{\lambda}_i/\lambda_i$ in Theorem~\ref{th1.2} does not mean that $(1-y)\hat \lambda_i$ is a good estimator of $\lambda_i$. In this section, we establish the CLT for $\hat\lambda_i$ to provide further properties.

We first introduce a centered parameter for $\hat\lambda_i$.
Let $\theta_{i}\in \mathbb{R}$, $1\leq i \leq q$, satisfy
\begin{align}\label{def_theta}
 1-\frac{1}{n}{\rm E}\left[\text{tr}\left\{\mathbf{M}^{-1}(\theta_{i})\right\}\right]=\frac{\lambda_{i}}{\theta_{i}}\left( 1+\frac{1}{T}{\rm E}\left[\text{tr}\left\{\mathbf{M}^{-1}\left(\theta_{i}\right)\frac{\mathbf{F}_{0}}{\theta_{i}}\right\}\right]\right),
\end{align}
and define $\delta_i$, for  $1\le i\le q$, as
\begin{align}\label{100}
\delta_{i}=\frac{\hat\lambda_{i}-\theta_{i}}{\theta_{i}}.
\end{align}
By Lemma~\ref{fisherlemma2}, when $n\rightarrow\infty$, we can easily see that
\begin{align*}
\frac{1}{p-q}{\rm E} \left[\text{tr} \left\{\mathbf{M}^{-1}(\theta_{i})\right\}\right]={\rm E}\left\{\widetilde{ m}_{\theta_{i}}(1)\right\}\rightarrow 1\quad \text{and}\quad
\frac{1}{p-q}{\rm E}\left[\text{tr}\left\{\mathbf{M}^{-1}(\theta_{i})\frac{\mathbf{F}_{0}}{\theta_{i}}\right\}\right]\rightarrow 0.
\end{align*}
It follows by (\ref{def_theta}) that
\begin{align*}
\frac{\lambda_{i}}{\theta_{i}}=\left(1-\frac{p-q}{n}\right)+{\rm o}(1)\rightarrow 1-y.
\end{align*}
Since the equation in  Definition~\ref{def_theta} for $\theta_i$ is  hard to calculate, an alternative definition for $\theta_i$ is proposed as follows.
Recall the definition of  $\widetilde m_{\theta}(z)$ in (\ref{def_m}):
\begin{align*}
\widetilde{ m}_{\theta}(z)&=\frac{1}{p-q}\text{tr} \left(z\textbf{I}_{p-q}-\frac{\textbf{F}_{0}}{\theta}  \right)^{-1},\ \theta\in \mathbb{R},\ z\in \mathbb{C}^{+}.
\end{align*}
Denoting $f_{\theta}(x)=\theta/(\theta-x)$ for any fixed $\theta\in \mathbb{R}$,  we have
\begin{gather*}
\widetilde{ m}_{\theta}(1)=\frac{1}{p-q}  \text{tr} \left(\textbf{I}_{p-q}-\frac{\textbf{F}_{0}}{\theta}  \right)^{-1}=\int_{-\infty}^{\infty}\frac{\theta}{\theta-x}dF_n(x)
=:F_n(f_{\theta}),
\end{gather*}
where $F_n$ denotes the ESD of the matrix $\mathbf{F}_0$.
By the CLT for linear spectral statistics (LSS) of Fisher matrices (see Theorem 3.10 in \cite{yao2015large}), for any fixed $\theta$,
\begin{gather*}
p\{F_n(f_{\theta})-F_{\tilde{c}_p,\tilde{y}_p}(f_{\theta})\}
\end{gather*}
converges weakly to a Gaussian variable. 
It follows that
\begin{align*}
&\widetilde{ m}_{\theta}(1)=F_{\tilde{c}_p,\tilde{y}_p}(f_{\theta})+{\rm O}_p(n^{-1})=-\theta \widetilde{\mathcal{S}}(\theta)+{\rm O}_p(n^{-1})\\
=&\frac{\tilde{c}_p-1}{\tilde{c}_p}+\frac{ \tilde{c}_p\{\theta(1-\tilde{y}_p)+1-\tilde{c}_p\}+2\theta \tilde{y}_p-\tilde{c}_p\sqrt{\{\theta(1-\tilde{y}_p)+1-\tilde{c}_p\}^{2}-4\theta}}{2\tilde{c}_p(\tilde{c}_p+\theta \tilde{y}_p)}+{\rm O}_p(n^{-1}),
\end{align*}
where $\widetilde{\mathcal{S}}(\cdot)$ denotes the stieltjes transform of $F_{\tilde{c}_p,\tilde{y}_p}$. This leads to
\begin{align}\label{ms}
{\rm E}\left\{\widetilde{ m}_{\theta}\left(1\right)\right\}=-\theta \widetilde{\mathcal{S}}(\theta)+{\rm O } (n^{-1}).
\end{align}
The definition of $\theta_{i}$ in (\ref{def_theta}) can be rewritten as
\begin{align*}
 1-\tilde{y}_p{\rm E}\left\{\widetilde{m}_{\theta_{i}}(1)\right\}=\frac{\lambda_{i}}{\theta_{i}}\left[ 1-\tilde{c}_p+\tilde{c}_p{\rm E}\left\{\widetilde{m}_{\theta_{i}}(1)\right\}\right].
\end{align*}
According to (\ref{ms}), it is equivalent to
\begin{align}\label{def_theta20}
 1+\tilde{y}_p\theta_{i}\widetilde{\mathcal{S}}(\theta_{i})+{\rm O}(n^{-1})=\frac{\lambda_{i}}{\theta_{i}}\left\{ 1-\tilde{c}_p-\tilde{c}_p\theta_{i}\widetilde{\mathcal{S}}(\theta_{i})+{\rm O}(n^{-1})\right\}.
\end{align}
Thus, we give another definition of $\theta_i$ by the following equation
\begin{align}\label{def_theta2}
 1+\tilde{y}_p\theta_{i}\widetilde{\mathcal{S}}(\theta_{i})=\frac{\lambda_{i}}{\theta_{i}}\left\{ 1-\tilde{c}_p-\tilde{c}_p\theta_{i}\widetilde{\mathcal{S}}(\theta_{i})\right\}.
\end{align}
It is notable that the $\theta_i$ defined by (\ref{def_theta2}) is also applicable to the CLT of $\delta_i$ in the later section. Comparing two equations (\ref{def_theta20}) and  (\ref{def_theta2}), we can derive that the difference between two $\delta_i$'s respectively derived from these two equations is at most ${\rm O}(n^{-1})$, which is smaller than the scale $n^{-1/2}$ of $\delta_i$.
Even Taylor's expansion on the stieltjes transformantion  $\widetilde{\mathcal{S}}(\cdot)$ can be simply used to the equation (\ref{def_theta2}) and then get  the explicit forms of $\theta_i$, although some errors would appear. In the remaining parts of this paper, we use $\theta_i$ defined by (\ref{def_theta}) in all results and their proofs.

Consider the case where all the spiked eigenvalues are simple, that is, $n_{i}=1$ for all $1\leq i \leq \ell$, which means that $\Lambda_{1}={\rm diag}(\lambda_{1}, \lambda_{2}, \ldots, \lambda_{q})$.

\begin{theorem}\label{th1.3}
Under Assumptions~\ref{A1}, \ref{A2}, \ref{A3}, \ref{A4} and that $n_i=1,\ 1\le i\le \ell $, i.e., $\ell=q$, it holds that, for all $1\leq i \leq q$,
\begin{align*}
\sqrt{p}\frac{\delta_{i}}{\sigma_i} \xrightarrow{d} \mathcal N\left(0,1\right)
\end{align*}
with $\sigma_i^2:=(y+c)\nu_i-c-y(1-3y)(1-y)^{-1}$,
where
$\nu_{i}={\rm E}|\mathbf{u}_i^{\top}\mathbf{Z}_1\mathbf{e}_1|^{4}$, $\mathbf{e}_1=(1,0,\ldots,0)^\top\in \mathbb{R}^q$ and $\mathbf{u}_i\in \mathbb{R}^q$ is the $i$-th column of the matrix $\mathbf{U}^{\top}$.
\end{theorem}
\begin{remark}
When  the value of the variance $\sigma_i^2$ at the population level is unknown, for statistical inference, estimating $\sigma_i^2$ is in need. A natural estimation way would be to estimate the eigenvector $\mathbf{u}_i$ first.
For the spiked population model, \cite{cai2017limiting} shows that when a leading eigenvalue of the sample covariance matrix is divergent, the corresponding sample eigenvector is a good estimator for its population counterpart in terms of their inner product. However, the situation becomes much more difficult when it comes to the spiked Fisher matrix. Recalling the assumed structure $\Sigma_2^{-1/2}\Sigma_1\Sigma_2^{-1/2}=\mathbf{I}_p+\Delta$, we suppose that  $\mathbf{v}_i:=(\mathbf{u}_i^\top, 0,\ldots,0)^\top\in\mathbb{R}^p$ is the eigenvector of $\Sigma_2^{-1/2}\Sigma_1\Sigma_2^{-1/2}=\mathbf{I}_p+\Delta$ corresponding to $\lambda_i$ and $\hat{\mathbf{v}}_i$ is that of $\mathbf{S}_1=\Sigma_1^{1/2}\mathbf{Y}\mathbf{Y}^\top\Sigma_1^{1/2}$. Then $\Sigma_2^{1/2}\hat{\mathbf{v}}_i$ is the eigenvector of $(\mathbf{I}_p+\Delta)^{1/2}\mathbf{Y}\mathbf{Y}^\top(\mathbf{I}_p+\Delta)^{1/2}$ corresponding to the $i$-th largest eigenvalue. If $\Sigma_2$ is known or can be consistently estimated, $\Sigma_2^{1/2}\hat{\mathbf{v}}_i$ is a good estimator of $\mathbf{v}_i$, by Theorem~4.1 in \cite{cai2017limiting}. But actually $\Sigma_2$ cannot be easily recovered based on $\mathbf{S}_2$ because of the delocalization of those eigenvectors for non-outliers (see \cite{bloemendal2016principal}). Thus, how to construct a consistent estimation of $\Sigma_2$ becomes a challenging issue.  As a special case, when entries of $\mathbf{Y}$ and $\mathbf{Z}$ are Gaussian, the parameter $\nu_i$ equals to $3$, which is independent of the value of $\mathbf{u}_i$.
In practice, the bootstrap approximation would be an alternative way to achieve a reliable estimation of $\sigma_i^2$.
For estimation of the variance of the largest sample eigenvalue in a spiked population model, spiked population model, \cite{el2019non} shows that the bootstrap approximation works when the largest eigenvalue is quite large. This deserves a further study.
\end{remark}

\begin{figure}[h]
    \centering
     \begin{subfigure}
         \centering
         \includegraphics[scale=0.45]{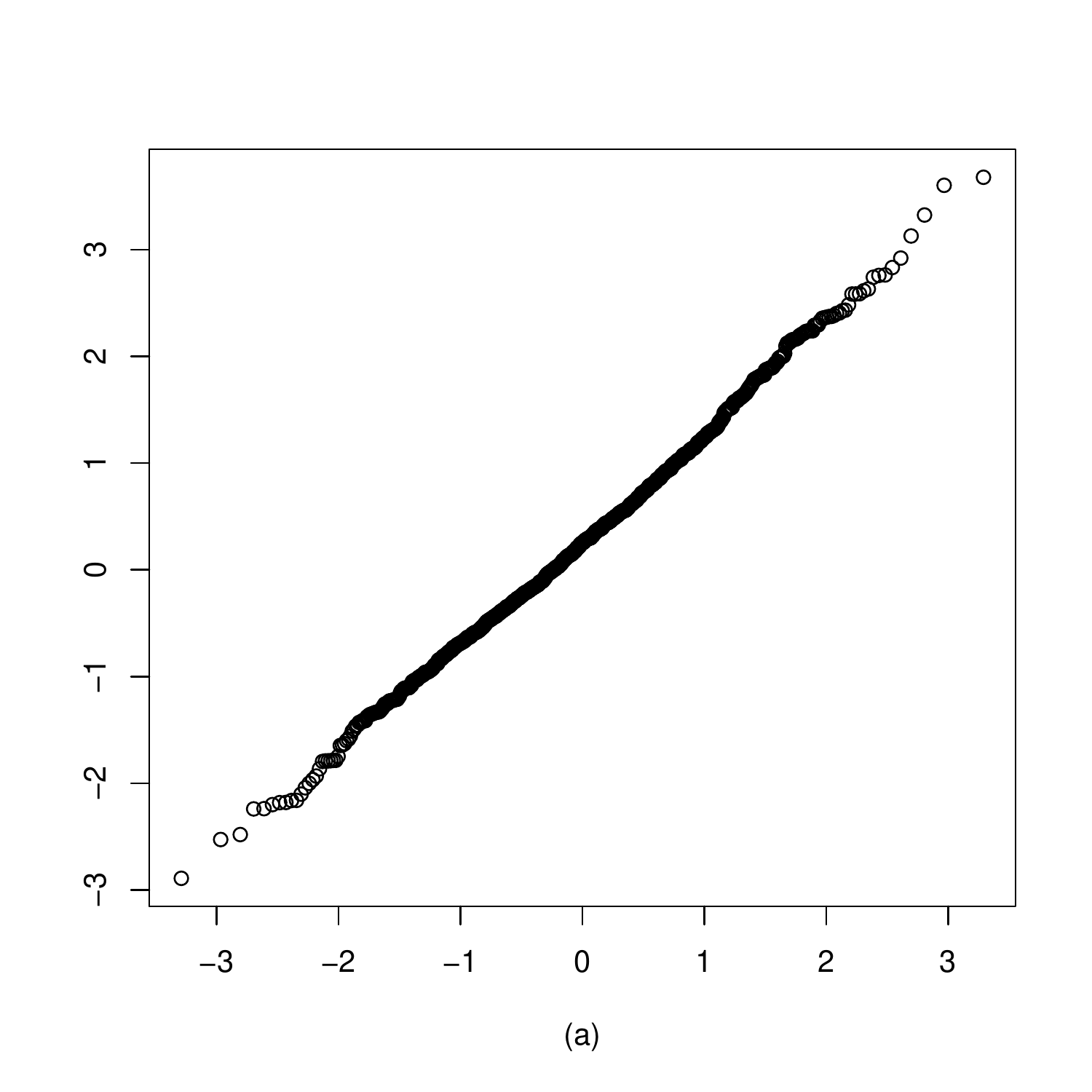}
     \end{subfigure}
     \hfill
     \begin{subfigure}
         \centering
         \includegraphics[scale=0.45]{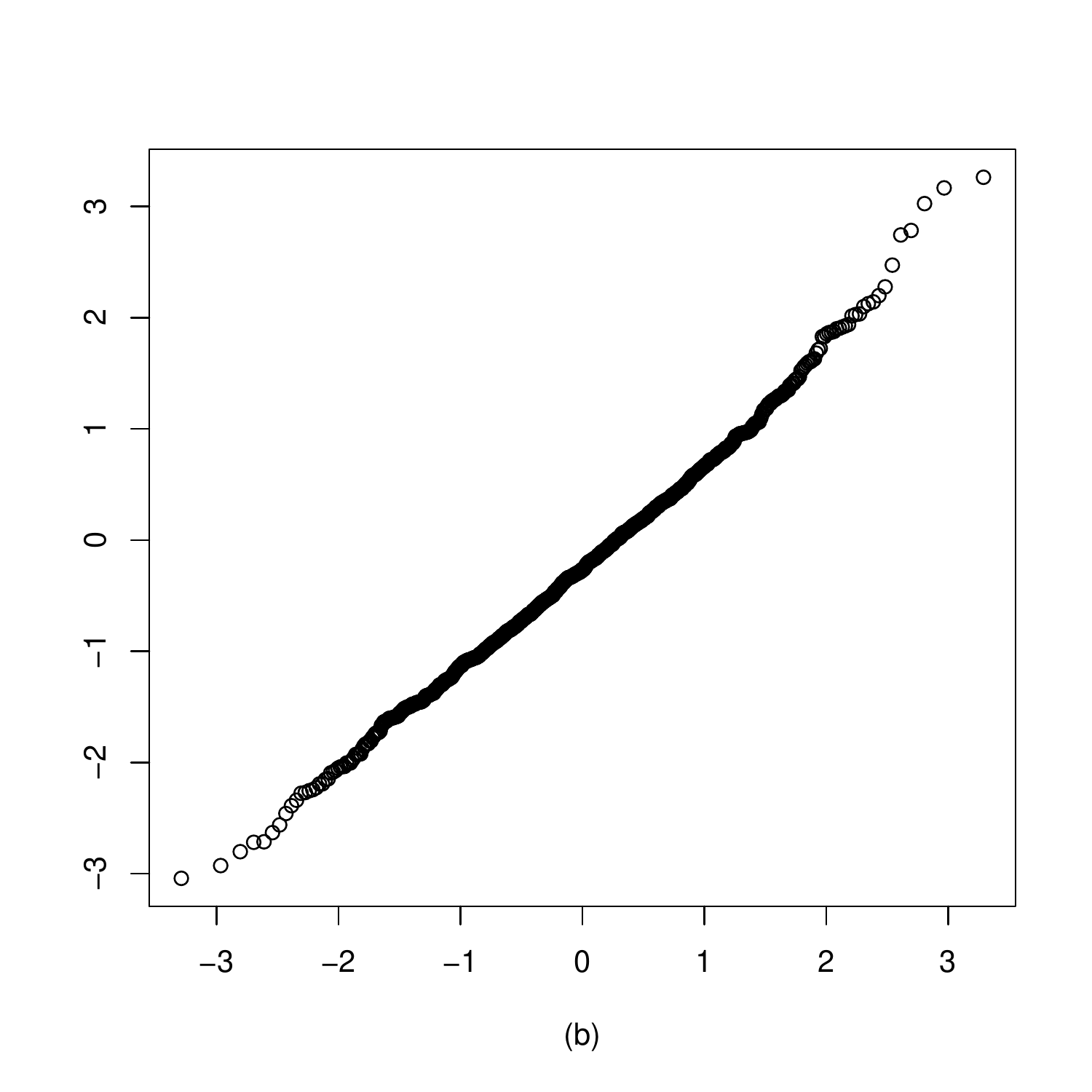}
     \end{subfigure}
     \caption{(a) The qq plot of the normalized largest spiked eigenvalue $\sqrt{p}\delta_1/\sigma_1$ from $1000$ independent replications. (b) The qq plot of the normalized smallest spiked eigenvalue $\sqrt{p}\delta_q/\sigma_q$ from $1000$ independent replications.}\label{qqplot}
\end{figure}

To check the practical applicability of Theorem~\ref{th1.3}, a simulation is conducted.
Set $p=200$, $T=600$, $n=1000$, 
$q=\lceil 2\log p  \rceil$, $\lambda_i=(3/2)^{q+1-i}(\log p/3)^3$ for $1\le i\le q$, where $\lceil x \rceil$ denotes the smallest integer greater than or equal to $x$. Let $\Sigma_1={\rm diag} (\lambda_1,\ldots,\lambda_q,1,\ldots,1)$ and $\Sigma_2=\mathbf{I}_p$.
Draw a sample $\{\mathbf x_i\}_{1\le i\le T}$ of size $T$ from $\mathcal{N}(0, \Sigma_1)$ and a sample $\{\mathbf z_i\}_{1\le i\le n}$ of size $n$ from $\mathcal{N}(0, \Sigma_2)$. Compute the largest $q$ eigenvalues $\hat{\lambda}_i$, $1\le i\le q$, of the Fisher matrix $\mathbf{F}=\mathbf{S}_2^{-1}\mathbf{S}_1$ and then $\delta_i$ accordingly, where $\mathbf{S}_1=\sum_{i=1}^T\mathbf{x}_i\mathbf{x}_i^\top/T$ and $\mathbf{S}_2=\sum_{i=1}^n\mathbf{z}_i\mathbf{z}_i^\top/n$. We draw qq plots of $\sqrt{p}\delta_1/\sigma_1$ and $\sqrt{p}\delta_q/\sigma_q$ from $1000$ independent replications in Figure~\ref{qqplot}. It suggests that both of $\sqrt{p}\delta_1/\sigma_1$ and $\sqrt{p}\delta_q/\sigma_q$ are well approximated by the standard normal distribution.\\

Next, consider the case where  some spiked eigenvalues are possibly  multiple:
\begin{align*}
\Lambda_{1}=\text{diag}(\underbrace{\lambda_{1},\ldots,\lambda_{N_1}}_{n_{1}},\ldots,\underbrace{\lambda_{N_{\ell-1}+1},\ldots,\lambda_{q}}_{n_{\ell}}),
\end{align*}
where $\lambda_1=\ldots=\lambda_{N_1}>\ldots>\lambda_{N_{\ell-1}+1}=\ldots=\lambda_{q}$, $N_i:=\sum_{j=1}^i n_j$ for $1\le i\le \ell$ and there exists a constant $C< \infty$ such that $1\leq n_{i}\leq C$ for all $1\le i\le \ell$.
According to the multiplicities of spiked eigenvalues, we divide  the index set $\{1,\ldots,q\}$ into $\ell$ subsets, $J_{i}=\left\{N_{i-1}+1,\ldots,N_{i}\right\},\ 1\le i\le \ell$. Here we denote $N_0=0$.
For any $1\le i\le \ell$, and $1\le h,k, h_1,k_1,h_2,k_2\le n_i$,  define
\begin{gather*}
\mathcal{M}_{N_i,h,k}:={\rm E}\left(\mathbf{u}_{N_{i-1}+h}^{\top}\mathbf{Z}_1\mathbf{e}_1\mathbf{u}_{N_{i-1}+k}^{\top}\mathbf{Z}_1\mathbf{e}_1\right),\\\mathcal{M}_{N_i,h_1,k_1,h_2,k_2}:={\rm E}\left(\mathbf{u}_{N_{i-1}+h_1}^{\top}\mathbf{Z}_1\mathbf{e}_1\mathbf{u}_{N_{i-1}+k_1}^{\top}\mathbf{Z}_1\mathbf{e}_1\mathbf{u}_{N_{i-1}+h_2}^{\top}\mathbf{Z}_1\mathbf{e}_1\mathbf{u}_{N_{i-1}+k_2}^{\top}\mathbf{Z}_1\mathbf{e}_1\right).
\end{gather*}
\begin{theorem}\label{th4.1}
Suppose that Assumptions  \ref{A1}, \ref{A2}, \ref{A3}, \ref{A4} and \ref{A5} hold.
Define  $\phi_{i}(\hat\lambda_{j})=(\hat\lambda_{j}-\theta_{j})/\theta_{j}$, for $1\le i\le \ell$ and $j\in J_{i}$. Then $\sqrt{p}\{\phi_{i}(\hat\lambda_{j}), j\in J_{i}    \}$ converges weakly to the distribution of the eigenvalues of the $n_{i}\times n_{i}$ random matrix $\Re^{(i)}$, where
$\Re^{(i)}=\left(R_{hk}^{(i)}\right)_{1\le h,k\le n_i}$ is a symmetric matrix with independent Gaussian entries of mean zero and covariance structure
\begin{align*}
{\rm cov}\left(R_{h_1,k_1}^{(i)},R_{h_2,k_2}^{(i)}\right)=&(1-y)^{-2}\omega\left(\mathcal{M}_{N_{i},h_1,k_1,h_2,k_2}-\mathcal{M}_{N_{i},h_1,k_1}\mathcal{M}_{N_{i},h_2,k_2}\right)\\
&+(1-y)^{-2}\left(\beta-\omega\right)\left(\mathcal{M}_{N_{i},h_1,k_2}\mathcal{M}_{N_{i},h_2,k_1}+\mathcal{M}_{N_{i},h_1,h_2}\mathcal{M}_{N_{i},k_1,k_2}\right),
\end{align*}
where $\omega =\left(y+c\right)\left(1-y\right)^2$ and $\beta = y\left(1-y\right)+c\left(1-y\right)^2$.
\end{theorem}

\section{Proofs of the theorems}\label{sec3}
We begin with a summary of the proofs.
Roughly, the proof of Theorem~\ref{th1.2} proceeds in three steps. First, we prove that the spiked eigenvalue $\hat{\lambda}_i$, $1\le i\le q$, solves the equation \eqref{object9091} whose left-hand side is the determinant of a $q\times q$ matrix which can be decomposed into four terms, namely $\mathbf{U}\Xi_A\mathbf{U}^\top$, $\mathbf{U}\Xi_B\mathbf{U}^\top$, $\mathbf{U}\Xi_C\mathbf{U}^\top$ and $\mathbf{U}\Xi_D\mathbf{U}^\top$ defined below. Second, we derive the limit of each entry of these four matrices and  their convergence rates in $\ell_\infty$ norm, where the CLT for random sequilinear forms in \cite{bai2008central} and Chebyshev's inequality are repeatedly used. Third, using eigenvalue perturbation theorems on \eqref{object9091}, we estimate the fluctuation of the scaled eigenvalue $\hat{\lambda}_i/\lambda_i$ and reach  the result. As for the proof of Theorem~\ref{th1.3}, we also work on the equation \eqref{object9091} in three main steps. First, we rewrite the matrix in \eqref{object9091} as the sum of $\mathbf{U}\Theta_{1n}\mathbf{U}^\top$, $\mathbf{U}\delta_i\Theta_{2n}\mathbf{U}^\top$ and $\mathbf{U}\Theta_{3n}\mathbf{U}^\top$. See equation~\eqref{object2} below. Second, we prove the CLT for each diagonal entry of $\mathbf{U}\Theta_{1n}\mathbf{U}^\top$ (Lemma~\ref{lemma1.1}) and estimate the $\ell_\infty$ norm of $\mathbf{U}\Theta_{1n}\mathbf{U}^\top$ (Lemma~\ref{lemma1.10}), $\mathbf{U}\Theta_{2n}\mathbf{U}^\top$ (Lemma~\ref{fisherlemma3}) and $\mathbf{U}\Theta_{3n}\mathbf{U}^\top$. Third, we expand the determinant in \eqref{object2} by Leibniz formula and then achieve the CLT for $\delta_i$. In this section, we will cite the lemmas given in the next section without the proofs whose details are postponed to the next section.\\

\noindent{\bf {Proof of Theorem \ref{th1.2}.}}
We first show that for $1\le i\le q$, $\hat\lambda_i$ converges to infinity at the same order with $\lambda_i$ almost surely, i.e., there exists some constant $C>1$ such that $C^{-1}< \hat \lambda_i/\lambda_i<C$ almost surely.

For any $1\le i\le q$, by Theorem~A.10 in \cite{bai2010spectral}, we have that
\begin{align*}
\hat\lambda_i= s_i(\mathbf{S}_2^{-1}\mathbf{S}_1)\le s_i(\mathbf{S}_1) s_1(\mathbf{S}_2^{-1})=s_i(\mathbf{S}_1)s_p^{-1}(\mathbf{S}_2)\quad \text{and}\quad s_i(\mathbf{S}_1)\le s_i(\mathbf{S}_2^{-1}\mathbf{S}_1) s_1\left(\mathbf{S}_2\right).
\end{align*}
Noting a basic fact that $s_1\left(\mathbf{S}_2\right)\rightarrow (1+\sqrt y)^2$ and $s_p\left(\mathbf{S}_2\right)\rightarrow (1-\sqrt y)^2>0$ almost surely, we have $0< C_1< \hat\lambda_i/s_i(\mathbf{S}_1)\le C_2< +\infty$ almost surely for some constants $C_1$ and $C_2$.

Again, by Theorem~A.10 in \cite{bai2010spectral} and Weyl's inequality, we have
\begin{align*}
s_i(\mathbf{S}_1)\le s_i(\Sigma_1)s_1\left(\frac{1}{T}\mathbf{Y}\mathbf{Y}^\top\right)=\lambda_i s_1\left(\frac{1}{T}\mathbf{Y}\mathbf{Y}^\top\right)
\end{align*}
and
\begin{align*}
s_i(\mathbf{S}_1)&=s_i\left(\frac{1}{T}\mathbf{Y}^\top\Sigma_1\mathbf{Y}\right)=s_i\left(\frac{1}{T}\mathbf{Y}_1^\top\Sigma_{11}\mathbf{Y}_1+\frac{1}{T}\mathbf{Y}_2^\top\mathbf{Y}_2\right)\geq s_i\left(\frac{1}{T}\mathbf{Y}_1^\top\Sigma_{11}\mathbf{Y}_1\right)\\
&\geq s_i\left(\Sigma_{11}\right)s_q\left(\frac{1}{T}\mathbf{Y}_1\mathbf{Y}_1^\top\right)=\lambda_i s_q\left(\frac{1}{T}\mathbf{Y}_1\mathbf{Y}_1^\top\right).
\end{align*}
Due to the fact that
\begin{align*}
s_1\left(\frac{1}{T}\mathbf{Y}\mathbf{Y}^\top\right)\rightarrow (1+\sqrt c)^2\quad \text{and}\quad s_q\left(\frac{1}{T}\mathbf{Y}_1\mathbf{Y}_1^\top\right)\rightarrow 1
\end{align*}
almost surely, we have $0<C_3<s_i(\mathbf{S}_1)/\lambda_i<C_4<+\infty$ almost surely for some constants $C_3$ and $C_4$.

Thus, we conclude that  $C^{-1}<\hat\lambda_i/\lambda_i<C$ almost surely for some constant $C$.

For any $1\leq i \leq q$, by the definition of $\hat\lambda_{i}$, it solves the equation $\det\left(\hat\lambda_{i}\textbf{I}-\textbf{S}_{2}^{-1}\textbf{S}_{1}\right)=0$,
or equivalently,
\begin{align}\label{eq1}
\det\left(\hat\lambda_{i}\textbf{S}_{2} - \textbf{S}_{1}\right)=0.
\end{align}
By the decomposition of $\textbf{S}_{1}$ and $\textbf{S}_{2}$ in \eqref{decomposition},  the equation (\ref{eq1}) can be rewritten as
\begin{align}\label{longeq}
\det \begin{pmatrix}\frac{\hat\lambda_{i}}{n}\mathbf{Z}_{1}\mathbf{Z}_{1}^{\top}-\frac{1}{T}\mathbf{X}_{1}\mathbf{X}_{1}^{\top}
&\frac{\hat\lambda_{i}}{n}\mathbf{Z}_{1}\mathbf{Z}_{2}^{\top}-\frac{1}{T}\mathbf{X}_{1}\mathbf{X}_{2}^{\top}
\\\frac{\hat\lambda_{i}}{n}\mathbf{Z}_{2}\mathbf{Z}_{1}^{\top}-\frac{1}{T}\mathbf{X}_{2}\mathbf{X}_{1}^{\top}
&\frac{\hat\lambda_{i}}{n}\mathbf{Z}_{2}\mathbf{Z}_{2}^{\top}-\frac{1}{T}\mathbf{X}_{2}\mathbf{X}_{2}^{\top}
\end{pmatrix}=0.
\end{align}
By the formula of the determinant of partitioned matrices, we know that $\det \begin{pmatrix}A
&B
\\C
&D
\end{pmatrix}=\det(D)\det(A-BD^{-1}C)$ when $D$ is nonsingular. As for $1\leq i \leq q$, $\hat\lambda_{i}$ is an outlier eigenvalue of $\textbf{S}_{2}^{-1}\textbf{S}_{1}$ because $\hat\lambda_{i}$ goes to infinity at the same order with $\lambda_i$, which means
\begin{align*}
\det\left(\frac{\hat\lambda_{i}}{n}\mathbf{Z}_{2}\mathbf{Z}_{2}^{\top}-\frac{1}{T}\mathbf{X}_{2}\mathbf{X}_{2}^{\top}\right)\neq 0,
\end{align*}
then it follows by (\ref{longeq}) that
{\small
\begin{align}\label{product}
\det\left\{\left(\frac{\hat\lambda_{i}}{n}\mathbf{Z}_{1}\mathbf{Z}_{1}^{\top}-\frac{1}{T}\mathbf{X}_{1}\mathbf{X}_{1}^{\top}\right)
-\left(\frac{\hat\lambda_{i}}{n}\mathbf{Z}_{1}\mathbf{Z}_{2}^{\top}-\frac{1}{T}\mathbf{X}_{1}\mathbf{X}_{2}^{\top}\right)
\left(\frac{\hat\lambda_{i}}{n}\mathbf{Z}_{2}\mathbf{Z}_{2}^{\top}-\frac{1}{T}\mathbf{X}_{2}\mathbf{X}_{2}^{\top}\right)^{-1}
\left(\frac{\hat\lambda_{i}}{n}\mathbf{Z}_{2}\mathbf{Z}_{1}^{\top}-\frac{1}{T}\mathbf{X}_{2}\mathbf{X}_{1}^{\top}\right)  \right\}=0.
\end{align}
}

For $\lambda\in \mathbb{R}$, defining
\begin{align*}
 &\mathbf{A}(\lambda)=\mathbf{Z}_{2}^{\top}\mathbf{M}^{-1}(\lambda)\left(\frac{1}{n}\mathbf{Z}_{2}\mathbf{Z}_{2}^{\top}\right)^{-1}\frac{1}{n}\mathbf{Z}_{2},\\
 &\mathbf{B}(\lambda)=\mathbf{X}_{2}^{\top}\mathbf{M}^{-1}(\lambda)\left(\frac{1}{n}\mathbf{Z}_{2}\mathbf{Z}_{2}^{\top}\right)^{-1}\frac{1}{\lambda T}\mathbf{X}_{2},\\
 &\mathbf{C}(\lambda)= \mathbf{Z}_{2}^{\top}\mathbf{M}^{-1}(\lambda)\left(\frac{1}{n}\mathbf{Z}_{2}\mathbf{Z}_{2}^{\top}\right)^{-1}\frac{1}{ \lambda T}\mathbf{X}_{2},\\
 &\mathbf{D}(\lambda)=\mathbf{X}_{2}^{\top}\mathbf{M}^{-1}(\lambda)\left(\frac{1}{n}\mathbf{Z}_{2}\mathbf{Z}_{2}^{\top}\right)^{-1}\frac{1}{\lambda n}\mathbf{Z}_{2},
\end{align*}
it holds that $\mathbf{A}\left(\lambda\right)=\mathbf{A}\left(\lambda\right)^{\top}$, $\mathbf{B}\left(\lambda\right)=\mathbf{B}\left(\lambda\right)^{\top}$ and $T\mathbf{C}\left(\lambda\right)=n\mathbf{D}\left(\lambda\right)^{\top}$.
Then some elementary calculations lead to
\begin{align}\label{object909}
\det\left[\hat\lambda_{i} \frac{\mathbf{{Z}}_{1}\left\{\textbf{I}_{n}-\mathbf{A}(\hat\lambda_{i})\right\}\mathbf{{Z}}_{1}^{\top}}{n}-\frac{\mathbf{{X}}_{1}\left\{\textbf{I}_{T}+\mathbf{B}(\hat\lambda_{i})\right\}\mathbf{{X}}_{1}^{\top}}{T}
+\hat\lambda_{i} \frac{\mathbf{{Z}}_{1}\mathbf{C}(\hat\lambda_{i})\mathbf{{X}}_{1}^{\top}}{n} +
\hat\lambda_{i} \frac{\mathbf{{X}}_{1}\mathbf{D}(\hat\lambda_{i})\mathbf{{Z}}_{1}^{\top}}{T}  \right]=0.
\end{align}
To ease the notation, we define
\begin{align*}
&\Xi_{A}:=\hat\lambda_{i} \frac{\mathbf{Z}_{1}\left\{\textbf{I}_{n}-\mathbf{A}(\hat\lambda_{i})\right\}\mathbf{Z}_{1}^{\top}}{n},\\
&\Xi_{B}:=\frac{\mathbf{X}_{1}\left\{\textbf{I}_{T}+\mathbf{B}(\hat\lambda_{i})\right\}\mathbf{{X}}_{1}^{\top}}{T},\\
&\Xi_{C}:=\hat\lambda_{i} \frac{\mathbf{Z}_{1}\mathbf{C}(\hat\lambda_{i})\mathbf{{X}}_{1}^{\top}}{n},\\
&\Xi_{D}:=\hat\lambda_{i} \frac{\mathbf{{X}}_{1}\mathbf{D}(\hat\lambda_{i})\mathbf{{Z}}_{1}^{\top}}{T}.
\end{align*}
Multiplying the matrix in (\ref{object909}) by $\mathbf{U}$ on the left side hand and by $\mathbf{U}^{\top}$ on the right side, we have
\begin{align}\label{object9091}
\det\left\{\mathbf{U}\left({\Xi_{A}}- {\Xi_{B}}+ {\Xi_{C}}+ {\Xi_{D}}\right)\mathbf{U}^{\top}\right\}=0.
\end{align}

Next, we  analyze these four terms in (\ref{object9091}) in the following.

For the term $\mathbf{U}{\Xi_{A}}\mathbf{U}^{\top}$,  we first consider the decomposition
\begin{align*}
\frac{1}{n}{\mathbf{{Z}}_{1}\left\{\textbf{I}_{n}-\mathbf{A}(\hat\lambda_{i})\right\}\mathbf{{Z}}_{1}^{\top}}=\frac{1}{n}{\mathbf{{Z}}_{1}\left\{\textbf{I}_{n}-\mathbf{A}(\lambda_{i})\right\}\mathbf{{Z}}_{1}^{\top}}+\frac{1}{n}{\mathbf{{Z}}_{1}\left\{\mathbf{A}(\hat\lambda_{i})-\mathbf{A}(\lambda_{i})\right\}\mathbf{{Z}}_{1}^{\top}}.
\end{align*}
By Lemma \ref{fisherlemma2} below,  we have
$\widetilde{ m}_{\lambda_{i}}(1) -1={\rm O}_{a.s.}(\lambda_{i}^{-1})$,
which implies
\begin{align*}
\frac{1}{n}\text{tr}\left\{\textbf{I}_{n}-\mathbf{A}(\lambda_{i})\right\}=1-\frac{p-q}{n}\widetilde{ m}_{ \lambda_{i}}(1)=1-y_p+\frac{q}{n}+{\rm O}_{a.s.}(\lambda_{i}^{-1}).
\end{align*}
Note that ${\rm E}(\mathbf{Z}_{1}\mathbf{Z}_{1}^{\top}/n)=\textbf{I}_{q}$ and that $(\mathbf{X}_1,\mathbf{Z}_1)$ is independent of $(\mathbf{X}_2,\mathbf{Z}_2)$. Under Assumption~\ref{A3}, by using Theorem 7.2 of \cite{bai2008central}, we have that, for all $1\leq j\leq q$,
\begin{align}\label{Adiagonal1}
\textbf{e}_{j}^{\top}\left[ \frac{1}{n}{\mathbf{{Z}}_{1}\left\{\textbf{I}_{n}-\mathbf{A}(\lambda_{i})\right\}\mathbf{{Z}}_{1}^{\top}}\right]\textbf{e}_{j}
-\left\{1-\frac{p-q}{n}\widetilde{ m}_{\lambda_{i}}(1)\right\}={\rm O}_p\left(\frac{1}{\sqrt n}\right)
\end{align}
and
\begin{align}\label{a1}
{\rm E}\left(\textbf{e}_{j}^{\top}\left[ \frac{1}{n}{\mathbf{{Z}}_{1}\left\{\textbf{I}_{n}-\mathbf{A}(\lambda_{i})\right\}\mathbf{{Z}}_{1}^{\top}}\right]\textbf{e}_{j}
-\left\{1-\frac{p-q}{n}\widetilde{ m}_{\lambda_{i}}(1)\right\}\right)^2={\rm O}\left(\frac{1}{n}\right)
\end{align}
for all $1\leq j\leq q$.
For those off-diagonal elements, we have that, for any $1\leq j_1\neq j_2\leq q$,
\begin{align}\label{a2}
\textbf{e}_{j_1}^{\top}\left[\frac{1}{n}\mathbf{{Z}}_1\left\{\textbf{I}_n-\mathbf{A}(\lambda_{i})\right\}\mathbf{{Z}}_1^{\top}\right]\textbf{e}_{j_2}
={\rm O}_{p}\left(\frac{1}{\sqrt n}\right)
\end{align}
and
\begin{align}\label{a3}
{\rm E}\left(\textbf{ e}_{j_1}^{\top}\left[\frac{1}{n}\mathbf{{Z}}_1\left\{\textbf{I}_n-\mathbf{A}(\lambda_{i})\right\}\mathbf{{Z}}_1^{\top}\right]\textbf{e}_{j_2}\right)^2
={\rm O}\left(\frac{1}{ n}\right),
\end{align}
which is implied by Theorem 7.1 and Corollary 7.1 in \cite{bai2008central}.
Also we can write
\begin{align*}
\mathbf{A}\left(\lambda_i\right)-\mathbf{A}\left(\hat\lambda_i\right)
&=\mathbf{Z}_{2}^{\top}\left\{\mathbf{M}^{-1}\left(\lambda_i\right)-\mathbf{M}^{-1}\left(\hat\lambda_i\right)\right\}\left(\frac{1}{n} \mathbf{Z}_{2} \mathbf{Z}_{2}^{\top}\right)^{-1} \frac{1}{n} \mathbf{Z}_{2}\nonumber\\
&=\mathbf{Z}_{2}^{\top}\mathbf{M}^{-1}\left(\lambda_i\right)\left\{\mathbf{M}\left(\hat\lambda_i\right)-\mathbf{M}\left(\lambda_i\right)\right\}\mathbf{M}^{-1}\left(\hat\lambda_i\right)\left(\frac{1}{n} \mathbf{Z}_{2} \mathbf{Z}_{2}^{\top}\right)^{-1} \frac{1}{n} \mathbf{Z}_{2}\nonumber\\
&= \left(\lambda_i^{-1}-\hat\lambda_i^{-1}\right)\mathbf{Z}_{2}^{\top}\mathbf{M}^{-1}\left(\lambda_i\right)\mathbf{F}_0\mathbf{M}^{-1}\left(\hat\lambda_i\right)\left(\frac{1}{n} \mathbf{Z}_{2} \mathbf{Z}_{2}^{\top}\right)^{-1} \frac{1}{n} \mathbf{Z}_{2}.
\end{align*}
It can be bounded by
\begin{align*}
&\left\|\mathbf{A}\left(\lambda_i\right)-\mathbf{A}\left(\hat\lambda_i\right)\right\|
= \left\|\left(\lambda_i^{-1}-\hat\lambda_i^{-1}\right)\mathbf{Z}_{2}^{\top}\mathbf{M}^{-1}\left(\lambda_i\right)\mathbf{F}_0\mathbf{M}^{-1}\left(\hat\lambda_i\right)\left(\frac{1}{n} \mathbf{Z}_{2} \mathbf{Z}_{2}^{\top}\right)^{-1} \frac{1}{n} \mathbf{Z}_{2}\right\|\nonumber\\
\le & \left|\lambda_i^{-1}-\hat\lambda_i^{-1}\right| \left\|\frac{1}{\sqrt n}\mathbf{Z}_{2}^{\top}\right\| \left\|\mathbf{M}^{-1}\left(\lambda_i\right)\right\| \left\|\mathbf{F}_0\right\| \left\|\mathbf{M}^{-1}\left(\hat\lambda_i\right)\right\| \left\|\left(\frac{1}{n} \mathbf{Z}_{2} \mathbf{Z}_{2}^{\top}\right)^{-1}\right\| \left\| \frac{1}{\sqrt n} \mathbf{Z}_{2}\right\|\nonumber
={\rm O}(\lambda_i^{-1})
\end{align*}
almost surely.
It follows that, for any $1\le j_1,j_2\le q$,
\begin{align}\label{Adifference0}
\mathbf{e}_{j_1}^{\top}\left[\frac{1}{n}{\mathbf{{Z}}_{1}\left\{\mathbf{A}(\hat\lambda_{i})-\mathbf{A}(\lambda_{i})\right\}\mathbf{{Z}}_{1}^{\top}}\right]\mathbf{e}_{j_2}={\rm O}_{a.s.}(\lambda_i^{-1}).
\end{align}

Combining \eqref{Adiagonal1}, \eqref{a2}) and \eqref{Adifference0}, we can get that, for any $1\le j\le q$,
\begin{align*}
\textbf{e}_{j}^{\top}\left[ \frac{1}{n}{\mathbf{{Z}}_{1}\left\{\textbf{I}_{n}-\mathbf{A}(\hat\lambda_{i})\right\}\mathbf{{Z}}_{1}^{\top}}\right]\textbf{e}_{j}
-\left(1-\frac{p-q}{n}\right)={\rm O}_p\left(\frac{1}{\sqrt n}\right)+{\rm O}_{a.s.}(\lambda_i^{-1})
\end{align*}
and that, for any $1\le j_1\neq j_2\le q$,
\begin{align*}
\textbf{ e}_{j_1}^{\top}\left[\frac{1}{n}\mathbf{{Z}}_1\left\{\textbf{I}_n-\mathbf{A}(\hat\lambda_{i})\right\}\mathbf{{Z}}_1^{\top}\right]\textbf{e}_{j_2}
={\rm O}_{p}\left(\frac{1}{\sqrt n}\right)+{\rm O}_{a.s.}(\lambda_i^{-1}).
\end{align*}

Replacing $\mathbf{{Z}}_1$ by $\mathbf{U}\mathbf{{Z}}_1$, it is easy to check that all the above conclusions  still hold:
\begin{align}\label{a4}
\mathbf{e}_{j}^{\top}\left[\frac{1}{n}{\mathbf{U}\mathbf{{Z}}_{1}\left\{\textbf{I}_{n}-\mathbf{A}(\hat\lambda_{i})\right\}\mathbf{{Z}}_{1}^{\top}\mathbf{U}^{\top}}\right]\mathbf{e}_{j} =1-\frac{p-q}{n}+{\rm O}_p\left(\frac{1}{\sqrt n}\right)+{\rm O}_{a.s.}(\lambda_i^{-1})
\end{align}
for all $1\leq j\leq q$, and
\begin{align}\label{a6}
\mathbf{e}_{j_1}^{\top}\left[\frac{1}{n}{\mathbf{U}\mathbf{{Z}}_{1}\left\{\textbf{I}_{n}-\mathbf{A}(\hat\lambda_{i})\right\}\mathbf{{Z}}_{1}^{\top}\mathbf{U}^{\top}}
\right]\mathbf{e}_{j_2}={\rm O}_{p}\left( \frac{1}{\sqrt{n}}\right)+{\rm O}_{a.s.}(\lambda_i^{-1})
\end{align}
for all $1\le j_1\neq j_2\le q$.
By the definition of $\Xi_{A}$ in (\ref{object909}), together with \eqref{a4} and \eqref{a6},  we can see that, for all $1\leq j\leq q$,
\begin{align}\label{a7}
\mathbf{e}_{j}^{\top}{\mathbf{U}\Xi_A \mathbf{U}^{\top}}\mathbf{e}_{j} =\hat\lambda_i\left( 1-\frac{p-q}{n}\right)+\lambda_i \cdot {\rm O}_p\left(\frac{1}{\sqrt n}\right)+{\rm O}_{a.s.}(1)
\end{align}
and that, for all $1\le j_1\neq j_2\le q$,
\begin{align}\label{a8}
\mathbf{e}_{j_1}^{\top}{\mathbf{U}\Xi_A \mathbf{U}^{\top}}\mathbf{e}_{j_2}=\lambda_i \cdot {\rm O}_p\left(\frac{1}{\sqrt n}\right)+{\rm O}_{a.s.}(1).
\end{align}

For the term  $ \mathbf{U}\Xi_{B}\mathbf{U}^{\top} $, by the definition of $\mathbf{X}_{1}$, we can derive that
\begin{align*}
\mathbf{U}\Xi_{B}\mathbf{U}^{\top}=&\frac{1}{T}\Lambda_{1}^{\frac{1}{2}}\mathbf{U}\mathbf{{Y}}_{1}\left\{\textbf{I}_{T}+\mathbf{B}(\hat\lambda_{i})\right\}\mathbf{{Y}}_{1}^{\top}\mathbf{U}^{\top}\Lambda_{1}^{\frac{1}{2}}= \frac{1}{T}\Lambda_{1}^{\frac{1}{2}}
\mathbf{U}\mathbf{{Y}}_{1}\left\{\textbf{I}_{T}+\mathbf{B}(\lambda_{i})\right\}\mathbf{{Y}}_{1}^{\top}\mathbf{U}^{\top}\Lambda_{1}^{\frac{1}{2}}\\
&+\frac{1}{T}\Lambda_{1}^{\frac{1}{2}}
\mathbf{U}\mathbf{{Y}}_{1}\left\{\mathbf{B}(\hat\lambda_{i})-\mathbf{B}(\lambda_{i})\right\}\mathbf{{Y}}_{1}^{\top}\mathbf{U}^{\top}\Lambda_{1}^{\frac{1}{2}},
\end{align*}
where
\begin{align*}
\frac{1}{T}\text{tr}\left\{\mathbf{I}_T+\mathbf{B}\left(\lambda_i\right)\right\}=&\frac{1}{T}\text{tr}\left\{\mathbf{I}_T+\mathbf{X}_{2}^{\top}\mathbf{M}^{-1}(\lambda_i)\left(\frac{1}{n}\mathbf{Z}_{2}\mathbf{Z}_{2}^{\top}\right)^{-1}\frac{1}{\lambda_i T}\mathbf{X}_{2}\right\}\\
=& 1+\frac{1}{T}\text{tr}\left\{\mathbf{M}^{-1}(\lambda_i)\left(\frac{1}{n}\mathbf{Z}_{2}\mathbf{Z}_{2}^{\top}\right)^{-1}\frac{1}{\lambda_i T}\mathbf{X}_{2}\mathbf{X}_{2}^{\top}\right\}\\
=&1+\frac{1}{T}\text{tr}\left\{\left(\mathbf{I}_{p-q}-\frac{\mathbf{F}_0}{\lambda_i}\right)^{-1}\frac{\mathbf{F}_0}{\lambda_i}\right\}\\
=&1+\frac{p-q}{T}\left\{\widetilde m_{\lambda_i}(1)-1\right\}
\end{align*}
and
\begin{align*}
\mathbf{B}(\hat\lambda_{i})-\mathbf{B}(\lambda_{i})=&\mathbf{X}_{2}^{\top}\left\{\hat\lambda_i^{-1}\mathbf{M}^{-1}\left(\hat\lambda_i\right)-\lambda_i^{-1}\mathbf{M}^{-1}(\lambda_i)\right\}\left(\frac{1}{n}\mathbf{Z}_{2}\mathbf{Z}_{2}^{\top}\right)^{-1}\frac{1}{ T}\mathbf{X}_{2}\\
=&\hat\lambda_i^{-1}\lambda_i^{-1}\mathbf{X}_{2}^{\top}\mathbf{M}^{-1}\left(\hat\lambda_i\right)\left\{\lambda_i\mathbf{M}\left(\lambda_i\right)-\hat\lambda_i\mathbf{M}\left(\hat\lambda_i\right)\right\}\mathbf{M}^{-1}(\lambda_i)\left(\frac{1}{n}\mathbf{Z}_{2}\mathbf{Z}_{2}^{\top}\right)^{-1}\frac{1}{ T}\mathbf{X}_{2}\\
=&\left(\hat\lambda_i^{-1}-\lambda_i^{-1}\right)\mathbf{X}_{2}^{\top}\mathbf{M}^{-1}\left(\hat\lambda_i\right)\mathbf{M}^{-1}(\lambda_i)\left(\frac{1}{n}\mathbf{Z}_{2}\mathbf{Z}_{2}^{\top}\right)^{-1}\frac{1}{ T}\mathbf{X}_{2}.
\end{align*}
The same arguments for deriving \eqref{a7} and \eqref{a8} lead to that, for all $1\leq j\leq q$,
\begin{align}
\mathbf{e}_{j}^{\top}{\mathbf{U}\Xi_B \mathbf{U}^{\top}}\mathbf{e}_{j} =\lambda_j+\lambda_j \cdot {\rm O}_p\left(\frac{1}{\sqrt n}\right)+\lambda_j \cdot {\rm O}_{a.s.}(\lambda_i^{-1})
\end{align}
and that, for $ 1\leq j_1,j_2\leq q$,
\begin{align}
\mathbf{e}_{j_1}^{\top}{\mathbf{U}\Xi_B \mathbf{U}^{\top}}\mathbf{e}_{j_2}=\lambda_{j_1}^{\frac{1}{2}}\lambda_{j_2}^{\frac{1}{2}} \cdot {\rm O}_p\left(\frac{1}{\sqrt n}\right)+\lambda_{j_1}^{\frac{1}{2}}\lambda_{j_2}^{\frac{1}{2}} \cdot {\rm O}_{a.s.}(\lambda_i^{-1})
\end{align}
for all $1\le j_1\neq j_2\le q$.

For the term $\mathbf{U}\left(\Xi_{C}+\Xi_{D}\right)\mathbf{U}^{\top}$, by using the fact that $\mathbf{Y}_{1}=\mathbf{U}^{\top}\Lambda_{1}^{-\frac{1}{2}}\mathbf{U}\mathbf{X}_{1}$, we have that
\begin{align*}
&\mathbf{U}\left(\Xi_{C}+\Xi_{D}\right)\mathbf{U}^{\top}\\
=&\mathbf{U}\left\{\hat\lambda_{i}\frac{\mathbf{Z}_1\mathbf{C}(\hat\lambda_{i})\mathbf{X}_1^{\top}}{n}+\hat\lambda_{i}\frac{\mathbf{X}_1\mathbf{D}(\hat\lambda_{i})\mathbf{Z}_1^{\top}}{T}\right\}\mathbf{U}^{\top}\\
=&\mathbf{U}\left\{\lambda_{i}\frac{\mathbf{Z}_1\mathbf{C}(\lambda_{i})\mathbf{X}_1^{\top}}{n}+\lambda_{i}\frac{\mathbf{X}_1\mathbf{D}(\lambda_{i})\mathbf{Z}_1^{\top}}{T}\right\}\mathbf{U}^{\top} + \mathbf{U}\mathbf{Z}_1\left\{\hat\lambda_{i}\mathbf{C}(\hat\lambda_{i})-\lambda_{i}\mathbf{C}\left(\lambda_{i}\right)\right\}\frac{1}{n}\mathbf{Y}_1^{\top}\mathbf{U}^{\top}\Lambda_{1}^{\frac{1}{2}}\\
&+\Lambda_{1}^{\frac{1}{2}}\mathbf{U}\mathbf{Y}_1\left\{\hat\lambda_{i}\mathbf{D}(\hat\lambda_{i})-\lambda_{i}\mathbf{D}\left(\lambda_{i}\right)\right\}\frac{1}{T}\mathbf{Z}_1^{\top}\mathbf{U}^{\top},
\end{align*}
and that
\begin{align*}
&\mathbf{U}\left\{\lambda_{i}\frac{\mathbf{Z}_1\mathbf{C}(\lambda_{i})\mathbf{X}_1^{\top}}{n}+\lambda_{i}\frac{\mathbf{X}_1\mathbf{D}(\lambda_{i})\mathbf{Z}_1^{\top}}{T}\right\}\mathbf{U}^{\top}\\
=&\mathbf{U}\begin{pmatrix} \mathbf{Z}_1&\mathbf{X}_1\end{pmatrix}\begin{pmatrix} O&\frac{\lambda_{i}\mathbf{C}(\lambda_{i})}{n}\\\frac{\lambda_{i}\mathbf{D}(\lambda_{i})}{T}&O \end{pmatrix}\begin{pmatrix} \mathbf{Z}_1^{\top}\\\mathbf{X}_1^{\top} \end{pmatrix}\mathbf{U}^{\top}\\
=&\begin{pmatrix} \mathbf{U}\mathbf{Z}_1&\Lambda_{1}^{\frac{1}{2}}\mathbf{U}\mathbf{Y}_1\end{pmatrix}\begin{pmatrix} O&\frac{\lambda_{i}\mathbf{C}(\lambda_{i})}{n}\\\frac{\lambda_{i}\mathbf{D}(\lambda_{i})}{T}&O \end{pmatrix}\begin{pmatrix} \mathbf{Z}_1^{\top}\mathbf{U}^{\top}\\ \mathbf{Y}_1^{\top} \mathbf{U}^{\top}\Lambda_{1}^{\frac{1}{2}}\end{pmatrix}.
\end{align*}
Then we have that, for all $1\leq j_1,j_2\leq q$,
\begin{align*}
&\mathbf{e}_{j_1}^{\top}\mathbf{U}\left\{\lambda_{i}\frac{\mathbf{Z}_1\mathbf{C}(\lambda_{i})\mathbf{X}_1^{\top}}{n}+\lambda_{i}\frac{\mathbf{X}_1\mathbf{D}(\lambda_{i})\mathbf{Z}_1^{\top}}{T}\right\}\mathbf{U}^{\top}\mathbf{e}_{j_2}\nonumber\\
=&\mathbf{e}_{j_1}^{\top}\begin{pmatrix} \mathbf{U}\mathbf{Z}_1&\Lambda_{1}^{\frac{1}{2}}\mathbf{U}\mathbf{Y}_1\end{pmatrix}\begin{pmatrix} O&\frac{\lambda_{i}\mathbf{C}(\lambda_{i})}{n}\\\frac{\lambda_{i}\mathbf{D}(\lambda_{i})}{T}&O \end{pmatrix}\begin{pmatrix} \mathbf{Z}_1^{\top}\mathbf{U}^{\top}\\ \mathbf{Y}_1^{\top} \mathbf{U}^{\top}\Lambda_{1}^{\frac{1}{2}}\end{pmatrix}\mathbf{e}_{j_2}\\
=&\mathbf{e}_{j_1}^{\top}\begin{pmatrix} \mathbf{U}\mathbf{Z}_1&\lambda_{j_1}^{\frac{1}{2}}\mathbf{U}\mathbf{Y}_1\end{pmatrix}\begin{pmatrix} O&\frac{\lambda_{i}\mathbf{C}(\lambda_{i})}{n}\\\frac{\lambda_{i}\mathbf{D}(\lambda_{i})}{T}&O \end{pmatrix}\begin{pmatrix} \mathbf{Z}_1^{\top}\mathbf{U}^{\top}\\ \lambda_{j_2}^{\frac{1}{2}}\mathbf{Y}_1^{\top} \mathbf{U}^{\top}\end{pmatrix}\mathbf{e}_{j_2}\\
=&\frac{\lambda_{j_1}^{\frac{1}{2}}+\lambda_{j_2}^{\frac{1}{2}}}{2}\cdot\mathbf{e}_{j_1}^{\top}\begin{pmatrix} \mathbf{U}\mathbf{Z}_1&\mathbf{U}\mathbf{Y}_1\end{pmatrix}\begin{pmatrix} O&\frac{\lambda_{i}\mathbf{C}(\lambda_{i})}{n}\\\frac{\lambda_{i}\mathbf{D}(\lambda_{i})}{T}&O \end{pmatrix}\begin{pmatrix} \mathbf{Z}_1^{\top}\mathbf{U}^{\top}\\ \mathbf{Y}_1^{\top} \mathbf{U}^{\top}\end{pmatrix}\mathbf{e}_{j_2}\\
&+\frac{\lambda_{j_1}^{\frac{1}{2}}-\lambda_{j_2}^{\frac{1}{2}}}{2\mathbf{i}}\cdot\mathbf{e}_{j_1}^{\top}\begin{pmatrix} \mathbf{U}\mathbf{Z}_1&\mathbf{U}\mathbf{Y}_1\end{pmatrix}\begin{pmatrix} O&\frac{\lambda_{i}\mathbf{C}(\lambda_{i})}{n}\cdot\mathbf{i}\\-\frac{\lambda_{i}\mathbf{D}(\lambda_{i})}{T}\cdot\mathbf{i}&O \end{pmatrix}\begin{pmatrix} \mathbf{Z}_1^{\top}\mathbf{U}^{\top}\\ \mathbf{Y}_1^{\top} \mathbf{U}^{\top}\end{pmatrix}\mathbf{e}_{j_2}\\
=&\frac{\lambda_{j_1}^{\frac{1}{2}}+\lambda_{j_2}^{\frac{1}{2}}}{2}\cdot{\rm O}_p\left(\frac{1}{\sqrt n}\right)+\frac{\lambda_{j_1}^{\frac{1}{2}}-\lambda_{j_2}^{\frac{1}{2}}}{2\mathbf{i}}\cdot {\rm O}_p\left(\frac{1}{\sqrt n}\right)
=\left(\lambda_{j_1}^{\frac{1}{2}}+\lambda_{j_2}^{\frac{1}{2}}\right)\cdot{\rm O}_p\left(\frac{1}{\sqrt n}\right),
\end{align*}
where $\mathbf{i}:=\sqrt {-1}$ is the imaginary unit and the penultimate equality is implied by Theorem~7.1 in \cite{bai2008central}. Due to the fact that
\begin{align*}
\hat\lambda_{i}\mathbf{C}(\hat\lambda_{i})-\lambda_{i}\mathbf{C}\left(\lambda_{i}\right)
=&\mathbf{Z}_{2}^{\top} \left\{\mathbf{M}^{-1}(\hat\lambda_i)-\mathbf{M}^{-1}(\lambda_i)\right\}\left(\frac{1}{n} \mathbf{Z}_{2} \mathbf{Z}_{2}^{\top}\right)^{-1} \frac{1}{ T} \mathbf{X}_{2}\\
=&\mathbf{Z}_{2}^{\top} \mathbf{M}^{-1}(\hat\lambda_i)\left\{\mathbf{M}(\lambda_i)-\mathbf{M}(\hat\lambda_i)\right\}\mathbf{M}^{-1}(\lambda_i)\left(\frac{1}{n} \mathbf{Z}_{2} \mathbf{Z}_{2}^{\top}\right)^{-1} \frac{1}{ T} \mathbf{X}_{2}\\
=&\left(\hat\lambda_i^{-1}-\lambda_i^{-1}\right)\mathbf{Z}_{2}^{\top} \mathbf{M}^{-1}(\hat\lambda_i)\mathbf{F}_0\mathbf{M}^{-1}(\lambda_i)\left(\frac{1}{n} \mathbf{Z}_{2} \mathbf{Z}_{2}^{\top}\right)^{-1} \frac{1}{ T} \mathbf{X}_{2},\\
\hat\lambda_{i}\mathbf{D}(\hat\lambda_{i})-\lambda_{i}\mathbf{D}\left(\lambda_{i}\right)
=&\mathbf{X}_{2}^{\top} \left\{\mathbf{M}^{-1}(\hat\lambda_i)-\mathbf{M}^{-1}(\lambda_i)\right\}\left(\frac{1}{n} \mathbf{Z}_{2} \mathbf{Z}_{2}^{\top}\right)^{-1} \frac{1}{n} \mathbf{Z}_{2}\\
=&\mathbf{X}_{2}^{\top} \mathbf{M}^{-1}(\hat\lambda_i)\left\{\mathbf{M}(\lambda_i)-\mathbf{M}(\hat\lambda_i)\right\}\mathbf{M}^{-1}(\lambda_i)\left(\frac{1}{n} \mathbf{Z}_{2} \mathbf{Z}_{2}^{\top}\right)^{-1} \frac{1}{n} \mathbf{Z}_{2}\\
=&\left(\hat\lambda_i^{-1}-\lambda_i^{-1}\right)\mathbf{X}_{2}^{\top} \mathbf{M}^{-1}(\hat\lambda_i)\mathbf{F}_0\mathbf{M}^{-1}(\lambda_i)\left(\frac{1}{n} \mathbf{Z}_{2} \mathbf{Z}_{2}^{\top}\right)^{-1} \frac{1}{n} \mathbf{Z}_{2},
\end{align*}
we can get
\begin{align*}
\mathbf{e}_{j_1}^{\top}\mathbf{U}\mathbf{Z}_1\left\{\hat\lambda_{i}\mathbf{C}(\hat\lambda_{i})-\lambda_{i}\mathbf{C}\left(\lambda_{i}\right)\right\}\frac{1}{n}\mathbf{Y}_1^{\top}\mathbf{U}^{\top}\Lambda_{1}^{\frac{1}{2}}\mathbf{e}_{j_2}&=\lambda_{j_2}^{\frac{1}{2}}\cdot {\rm O}_{a.s.}\left(\lambda_i^{-1}\right)
\end{align*}
and
\begin{align*}
\mathbf{e}_{j_1}^{\top}\Lambda_{1}^{\frac{1}{2}}\mathbf{U}\mathbf{Y}_1\left\{\hat\lambda_{i}\mathbf{D}(\hat\lambda_{i})-\lambda_{i}\mathbf{D}\left(\lambda_{i}\right)\right\}\frac{1}{T}\mathbf{Z}_1^{\top}\mathbf{U}^{\top}\mathbf{e}_{j_2}=\lambda_{j_1}^{\frac{1}{2}}\cdot {\rm O}_{a.s.}\left(\lambda_i^{-1}\right)
\end{align*}
for any $1\le j_1,j_2\le q$. By using the similar  arguments for proving \eqref{a7} and \eqref{a8}, it  holds that
\begin{align}\label{cd1}
\mathbf{e}_{j_1}^{\top}\mathbf{U}\left(\Xi_{C}+\Xi_{D}\right)\mathbf{U}^{\top}\mathbf{e}_{j_2}=\left(\lambda_{j_1}^{\frac{1}{2}}+\lambda_{j_2}^{\frac{1}{2}}\right)\cdot {\rm O}_p\left(\frac{1}{\sqrt n}\right)+\left(\lambda_{j_1}^{\frac{1}{2}}+\lambda_{j_2}^{\frac{1}{2}}\right)\cdot {\rm O}_{a.s.}\left(\lambda_i^{-1}\right)
\end{align}
for any $1\le j_1,j_2\le q$.

Combining \eqref{a8}-\eqref{cd1} and the determinant \eqref{object9091}, we can compute the limit of $\hat\lambda_{i}/\lambda_i$ for each $1\leq i\leq q$. We use a new notation to denote the matrix in the determinant (\ref{object9091}).
Define
\begin{align*}
\Xi:=\mathbf{U}\left({\Xi_{A}}- {\Xi_{B}}+ {\Xi_{C}}+ {\Xi_{D}}\right)\mathbf{U}^{\top},\quad
\widetilde\Xi:=\text{diag}\left(\xi_{11},\ldots,\xi_{qq}\right),
\end{align*}
where $\xi_{jj}=\hat\lambda_i\left\{ 1-\left(p-q\right)/n\right\}-\lambda_j$.
Then by \eqref{a8}-\eqref{cd1}, we have that
\begin{align*}
\mathbf{e}_{j_1}^{\top}\left(\Xi-\widetilde\Xi\right)\mathbf{e}_{j_2}=&\lambda_i \cdot {\rm O}_p\left(\frac{1}{\sqrt n}\right)+{\rm O}_{p}(1)+\lambda_{j_1}^{\frac{1}{2}}\lambda_{j_2}^{\frac{1}{2}} \cdot{\rm O}_p\left(\frac{1}{\sqrt n}\right)+\lambda_{j_1}^{\frac{1}{2}}\lambda_{j_2}^{\frac{1}{2}} \cdot {\rm O}_{a.s.}(\lambda_i^{-1})\\
&+\left(\lambda_{j_1}^{\frac{1}{2}}+\lambda_{j_2}^{\frac{1}{2}}\right)\cdot{\rm O}_p\left(\frac{1}{\sqrt n}\right)+\left(\lambda_{j_1}^{\frac{1}{2}}+\lambda_{j_2}^{\frac{1}{2}}\right)\cdot {\rm O}_{a.s.}\left(\lambda_i^{-1}\right)\\
=&\left(\lambda_i+\lambda_{j_1}^{\frac{1}{2}}\lambda_{j_2}^{\frac{1}{2}}\right)\left\{ {\rm O}_p\left(\frac{1}{\sqrt n}\right)+{\rm O}_{a.s.}\left(\lambda_i^{-1}\right)\right\}
\end{align*}
for any $1\le j_1,j_2\le q$, which follows that
\begin{align}\label{total}
\mathbf{e}_{j_1}^{\top}\lambda_i^{-1}\left(\Xi-\widetilde\Xi\right)\mathbf{e}_{j_2}
=&\left(1+\lambda_i^{-1}\lambda_{j_1}^{\frac{1}{2}}\lambda_{j_2}^{\frac{1}{2}}\right)\left\{ {\rm O}_p\left(\frac{1}{\sqrt n}\right)+{\rm O}_{a.s.}\left(\lambda_i^{-1}\right)\right\}.
\end{align}

According to \eqref{a1} and \eqref{a3} for $\Xi_{A}$ (similar results also hold for $\Xi_{B}$, $\Xi_{C}$ and $\Xi_{D}$), it can be easily checked that the variance of the term in \eqref{total} has the order
\begin{align*}
\left(1+\lambda_i^{-1}\lambda_{j_1}^{\frac{1}{2}}\lambda_{j_2}^{\frac{1}{2}}\right)^2 \left( n^{-\frac{1}{2}}+\lambda_i^{-1}\right)^2.
\end{align*}
By Chebyshev's inequality, we have that, for any $\epsilon>0$,
\begin{align*}
&\Pr \left\{\max_{1\leq j_{1}, j_{2} \leq q}\left|\mathbf{e}_{j_1}^{\top}\lambda_i^{-1}\left(\Xi-\widetilde\Xi\right)\mathbf{e}_{j_2}\right| \geq \epsilon \left(n^{-\frac{1}{2}}+\lambda_i^{-1}\right)\right\}\\
\leq &\sum_{1\leq j_{1}, j_{2}\leq q} \Pr \left\{ \left|\mathbf{e}_{j_1}^{\top}\lambda_i^{-1}\left(\Xi-\widetilde\Xi\right)\mathbf{e}_{j_2}\right|\geq \epsilon \left(n^{-\frac{1}{2}}+\lambda_i^{-1}\right)\right\}\nonumber\\
\leq& \sum_{1\leq j_{1}, j_{2}\leq q} \frac{{\rm E}\left\{\mathbf{e}_{j_1}^{\top}\lambda_i^{-1}\left(\Xi-\widetilde\Xi\right)\mathbf{e}_{j_2}\right\}^{2}}{\epsilon^{2}\left(n^{-\frac{1}{2}}+\lambda_i^{-1}\right)^2}\\
=&\sum_{1\leq j_{1}, j_{2}\leq q}\left(1+\lambda_i^{-1}\lambda_{j_1}^{\frac{1}{2}}\lambda_{j_2}^{\frac{1}{2}}\right)^2 \cdot {\rm O}\left(\epsilon^{-2}\right)
=\left(q+\lambda_i^{-1}\sum_{j=1}^q\lambda_j\right)^2\cdot {\rm O}\left(\epsilon^{-2}\right)
=\kappa_1^2 \cdot {\rm O}\left(\epsilon^{-2}\right),
\end{align*}
which means
\begin{align*}
\left\|\lambda_i^{-1}\left(\Xi-\widetilde\Xi\right)\right\|_\infty=\max_{1\leq j_{1}, j_{2} \leq q}\left|\mathbf{e}_{j_1}^{\top}\lambda_i^{-1}\left(\Xi-\widetilde\Xi\right)\mathbf{e}_{j_2}\right|
= \kappa_1\cdot {\rm O}_p\left(\frac{1}{\sqrt n}+\lambda_i^{-1}\right)
\end{align*}
and then
\begin{align*}
|||\lambda_i^{-1}(\Xi-\widetilde\Xi)|||_\infty\le q \|\lambda_i^{-1}(\Xi-\widetilde\Xi)\|_\infty
= \kappa_1  q \cdot {\rm O}_p\left(\frac{1}{\sqrt n}+\lambda_i^{-1}\right).
\end{align*}

Note that the determinant equation $\det \big(\widetilde\Xi\big)=0$ is equivalent to $\det \left(\lambda_i^{-1}\widetilde\Xi\right)=0$, that is,
\begin{align*}
\det \left\{\frac{\hat{\lambda}_i}{\lambda_i}\left(1-\frac{p-q}{n}\right)\mathbf{I}_q-\lambda_i^{-1} \Lambda_1\right\}=0.
\end{align*}
At the same time, the equation $\det \big(\Xi\big)=0$ is equivalent to $
\det \left(\lambda_i^{-1}\Xi\right)=0$, that is,
\begin{align*}
\det \left\{\frac{\hat{\lambda}_i}{\lambda_i}\left(1-\frac{p-q}{n}\right)\mathbf{I}_q-\lambda_i^{-1} \Lambda_1+\lambda_i^{-1}\left(\Xi-\widetilde\Xi\right)\right\}=0.
\end{align*}

%

By eigenvalue perturbation theorems (see Theorem~6.3.2 in Chapter~6, \cite{horn2012matrix}), we have
\begin{align*}
\left|\frac{\hat\lambda_{i}}{\lambda_i}\left(1-\frac{p-q}{n}\right)-1\right|\leq |||\lambda_i^{-1}(\Xi-\widetilde\Xi)|||_\infty=\kappa_1  q\cdot {\rm O}_p\left(\frac{1}{\sqrt n}+\lambda_i^{-1}\right),
\end{align*}
that is
\begin{align}\label{result1}
\frac{\hat\lambda_i}{\lambda_i}=\frac{1}{1-y}+{\rm O}\left(y_p-y\right)+\kappa_1  q\cdot {\rm O}_p\left(\frac{1}{\sqrt n}+\lambda_i^{-1}\right).
\end{align}

Instead, we can compare determinant equations
\begin{align*}
\det \left(\Lambda_1^{-\frac{1}{2}}\widetilde\Xi\Lambda_1^{-\frac{1}{2}}\right)=0\quad \textrm{and}\quad
\det \left(\Lambda_1^{-\frac{1}{2}}\Xi\Lambda_1^{-\frac{1}{2}}\right)=0,
\end{align*}
and then repeat all the derivations above to achieve an upper bound of $|||\Lambda_1^{-1/2}(\Xi-\widetilde\Xi)\Lambda_1^{-1/2}|||_\infty$. In this case, we can get
\begin{align}\label{result2}
\frac{\hat\lambda_i}{\lambda_i}=\frac{1}{1-y}+{\rm O}\left(y_p-y\right)+\kappa_2  q\cdot {\rm O}_p\left(\frac{1}{\sqrt n}+\lambda_i^{-1}\right).
\end{align}

Thus, \eqref{result1} and \eqref{result2} lead to
\begin{align*}
\frac{\hat\lambda_i}{\lambda_i}=\frac{1}{1-y}+{\rm O}\left(y_p-y\right)+\kappa  q\cdot {\rm O}_p\left(\frac{1}{\sqrt n}+\lambda_i^{-1}\right),
\end{align*}
where $\kappa q(n^{-1/2}+\lambda_i^{-1})={\rm o}(1)$ under Assumption~\ref{A2}.
The proof is finished.\quad $\square$

\vspace{0.2in}

\noindent{\bf {Proof of Theorem \ref{th1.3}.}}  We  begin with the equation on $\hat\lambda_{i}$ in (\ref{object909}).  Recall that we have expressed
(\ref{object909}) as
\begin{align}\label{object909909}
\det\big( \Xi_{A}-\Xi_{B}+\Xi_{C}+\Xi_{D}  \big)=0.
\end{align}

For the first term $\Xi_{A}$, we can write
\begin{align*}
\Xi_{A}=& \frac{\hat\lambda_{i}}{n} \mathbf{{Z}}_{1} \left[\left\{ \textbf{I}_{n}-\mathbf{A}(\hat\lambda_{i})\right\}-\left\{\textbf{I}_{n}-\mathbf{A}(\theta_{i})\right\}\right]\mathbf{{Z}}_{1}^{\top}  + \frac{\hat\lambda_{i}}{n} \left(\mathbf{{Z}}_{1}\left\{\textbf{I}_{n}-\mathbf{A}(\theta_{i})\right\}\mathbf{{Z}}_{1}^{\top}-{\rm E}\left[ \mathbf{{Z}}_{1}\left\{\textbf{I}_{n}-\mathbf{A}(\theta_{i})\right\}\mathbf{{Z}}_{1}^{\top} \right] \right) \\
&+  \frac{\hat\lambda_{i}}{n}{\rm E}\left[\mathbf{{Z}}_{1}\left\{\textbf{I}_{n}-\mathbf{A}(\theta_{i})\right\}\mathbf{{Z}}_{1}^{\top} \right].
\end{align*}
Using the fact
\begin{align*}
\left\{\textbf{I}_{n}-\mathbf{A}(\hat\lambda_{i})\right\}-\left\{\textbf{I}_{n}-\mathbf{A}(\theta_{i})\right\}=
-\delta_{i} \mathbf{A}(\theta_{i})+ \delta_{i}\mathbf{Z}_{2}^{\top}\mathbf{M}^{-1}(\hat\lambda_{i})\mathbf{M}^{-1}(\theta_{i})\left(\frac{1}{n}\mathbf{Z}_{2}\mathbf{Z}_{2}^{\top}\right)^{-1}\frac{1}{n}\mathbf{Z}_{2}
\end{align*}
and $\hat\lambda_{i}=\theta_{i}(1+\delta_{i})$ by (\ref{100}),
we can get
\begin{align}\label{seq111}
\Xi_{A}=&\theta_{i}\delta_{i}(1+\delta_{i}) \frac{1}{n} \left(\mathbf{{Z}}_{1}\left\{\textbf{I}_{n}-\mathbf{A}(\theta_{i})\right\}\mathbf{{Z}}_{1}^{\top}-{\rm E}\left[ \mathbf{{Z}}_{1}\left\{\textbf{I}_{n}-\mathbf{A}(\theta_{i})\right\}\mathbf{{Z}}_{1}^{\top} \right] \right) \nonumber\\
&+ \theta_{i}\delta_{i}(1+\delta_{i}) \frac{1}{n}\mathbf{{Z}}_{1}\mathbf{Z}_{2}^{\top}\mathbf{M}^{-1}(\hat\lambda_{i})\mathbf{M}^{-1}(\theta_{i})\left(\frac{1}{n}\mathbf{Z}_{2}\mathbf{Z}_{2}^{\top}\right)^{-1}\frac{1}{n}\mathbf{Z}_{2} \mathbf{{Z}}_{1}^{\top}    \nonumber\\
&-\theta_{i}\delta_{i}(1+\delta_{i}) \frac{1}{n} \mathbf{{Z}}_{1}\mathbf{{Z}}_{1}^{\top}+\theta_{i}\delta_{i}(1+\delta_{i}) \frac{1}{n} {\rm E}\left[ \mathbf{{Z}}_{1}\left\{\textbf{I}_{n}-\mathbf{A}(\theta_{i})\right\}\mathbf{{Z}}_{1}^{\top} \right]  \nonumber\\
&+ \theta_{i}(1+\delta_{i})\frac{1}{n} \left(\mathbf{{Z}}_{1}\left\{\textbf{I}_{n}-\mathbf{A}(\theta_{i})\right\}\mathbf{{Z}}_{1}^{\top}-{\rm E}\left[ \mathbf{{Z}}_{1}\left\{\textbf{I}_{n}-\mathbf{A}(\theta_{i})\right\}\mathbf{{Z}}_{1}^{\top} \right] \right)  \nonumber\\
&+ \theta_{i}(1+\delta_{i}) \frac{1}{n} {\rm E}\left[ \mathbf{{Z}}_{1}\left\{\textbf{I}_{n}-\mathbf{A}(\theta_{i})\right\}\mathbf{{Z}}_{1}^{\top} \right] \nonumber\\
=&  \theta_{i}(1+\delta_{i})^{2}\frac{1}{n} \left(\mathbf{{Z}}_{1}\left\{\textbf{I}_{n}-\mathbf{A}(\theta_{i})\right\}\mathbf{{Z}}_{1}^{\top}-{\rm E}\left[ \mathbf{{Z}}_{1}\left\{\textbf{I}_{n}-\mathbf{A}(\theta_{i})\right\}\mathbf{{Z}}_{1}^{\top} \right] \right)\nonumber\\
&+ \theta_{i}\delta_{i}(1+\delta_{i}) \frac{1}{n}\mathbf{{Z}}_{1}\mathbf{Z}_{2}^{\top}\mathbf{M}^{-1}(\hat\lambda_{i})\mathbf{M}^{-1}(\theta_{i})\left(\frac{1}{n}\mathbf{Z}_{2}\mathbf{Z}_{2}^{\top}\right)^{-1}\frac{1}{n}\mathbf{Z}_{2} \mathbf{{Z}}_{1}^{\top}  \nonumber\\
&-\theta_{i}\delta_{i}(1+\delta_{i}) \frac{1}{n} \mathbf{{Z}}_{1}\mathbf{{Z}}_{1}^{\top} +\theta_{i}(1+\delta_{i})^{2} \frac{1}{n} {\rm E}\left[ \mathbf{{Z}}_{1}\left\{\textbf{I}_{n}-\mathbf{A}(\theta_{i})\right\}\mathbf{{Z}}_{1}^{\top} \right] \nonumber\\
=:&\theta_{i}(1+\delta_{i})^{2}\Xi_{A1}+\theta_{i}\delta_{i}(1+\delta_{i})\Xi_{A2}-\theta_{i}\delta_{i}(1+\delta_{i})\Xi_{A3}
+\theta_{i}(1+\delta_{i})^{2}\Xi_{A4}.
\end{align}

For the second term  $\Xi_{B}$, we can similarly write
 \begin{align}\label{seq2}
\Xi_{B}=&\frac{1}{T}\mathbf{{X}}_{1}\left\{ \mathbf{B}(\hat\lambda_{i})-\mathbf{B}(\theta_{i})\right\}\mathbf{{X}}_{1}^{\top}+\frac{1}{T}\left( \mathbf{{X}}_{1} \left\{\textbf{I}_{T}+\mathbf{B}(\theta_{i})\right\}\mathbf{{X}}_{1}^{\top}-{\rm E}\left[\mathbf{{X}}_{1} \left\{\textbf{I}_{T}+\mathbf{B}(\theta_{i})\right\}\mathbf{{X}}_{1}^{\top} \right]\right)\nonumber\\
&+\frac{1}{T}{\rm E}\left[\mathbf{{X}}_{1} \left\{\textbf{I}_{T}+\mathbf{B}(\theta_{i})\right\}\mathbf{{X}}_{1}^{\top} \right]\nonumber\\
=&\frac{1}{T}\left( \mathbf{{X}}_{1} \left\{\textbf{I}_{T}+\mathbf{B}(\theta_{i})\right\}\mathbf{{X}}_{1}^{\top}-{\rm E}\left[\mathbf{{X}}_{1} \left\{\textbf{I}_{T}+\mathbf{B}(\theta_{i})\right\}\mathbf{{X}}_{1}^{\top} \right]\right)\nonumber\\
&-\frac{\delta_{i}}{\hat\lambda_{i}}\cdot \frac{1}{T}\mathbf{{X}}_{1}
\mathbf{X}_{2}^{\top}\mathbf{M}^{-1}(\hat\lambda_{i})\mathbf{M}^{-1}(\theta_{i})\left(\frac{1}{n}\mathbf{Z}_{2}\mathbf{Z}_{2}^{\top}\right)^{-1}\frac{1}{T}\mathbf{X}_{2}\mathbf{{X}}_{1}^{\top} +\frac{1}{T}{\rm E}\left[\mathbf{{X}}_{1} \left\{\textbf{I}_{T}+\mathbf{B}(\theta_{i})\right\}\mathbf{{X}}_{1}^{\top} \right]\nonumber\\
=:&\Xi_{B1}-\frac{\delta_{i}}{\hat\lambda_{i}}\Xi_{B2}+\Xi_{B3},
\end{align}
where the second equality above uses the fact
\begin{align*}
\mathbf{B}(\hat\lambda_{i})-\mathbf{B}(\theta_{i})=-\frac{\delta_{i}}{\hat\lambda_{i}}
\mathbf{X}_{2}^{\top}\mathbf{M}^{-1}(\hat\lambda_{i})\mathbf{M}^{-1}(\theta_{i})\left(\frac{1}{n}\mathbf{Z}_{2}\mathbf{Z}_{2}^{\top}\right)^{-1}\frac{1}{T}\mathbf{X}_{2}.
\end{align*}

For the term $\Xi_{C}$, we have
\begin{align*}
\Xi_{C}&= \frac{\hat\lambda_{i}}{n} \mathbf{{Z}}_{1} \left\{\mathbf{C}(\hat\lambda_{i})-\mathbf{C}(\theta_{i})\right\}\mathbf{{X}}_{1}^{\top}+  \frac{\hat\lambda_{i}-\theta_{i}}{n} \mathbf{{Z}}_{1} \mathbf{C}(\theta_{i})\mathbf{{X}}_{1}^{\top}+ \frac{\theta_{i}}{n}\left[ \mathbf{{Z}}_{1} \mathbf{C}(\theta_{i})\mathbf{{X}}_{1}^{\top}-{\rm E}\left\{\mathbf{{Z}}_{1} \mathbf{C}(\theta_{i})\mathbf{{X}}_{1}^{\top} \right\}\right].
\end{align*}
Using the fact
\begin{align*}
\mathbf{C}(\hat\lambda_{i})-\mathbf{C}(\theta_{i})=-\frac{\delta_{i}}{\hat\lambda_{i}}
\mathbf{Z}_{2}^{\top}\mathbf{M}^{-1}(\hat\lambda_{i})\mathbf{M}^{-1}(\theta_{i})\left(\frac{1}{n}\mathbf{Z}_{2}\mathbf{Z}_{2}^{\top}\right)^{-1}\frac{1}{T}\mathbf{X}_{2},
\end{align*}
we have the decomposition
\begin{align}\label{seq32}
\Xi_{C}&=\theta_{i}\frac{1}{n}\left[ \mathbf{{Z}}_{1} \mathbf{C}(\theta_{i})\mathbf{{X}}_{1}^{\top}-{\rm E}\left\{ \mathbf{{Z}}_{1} \mathbf{C}(\theta_{i})\mathbf{{X}}_{1}^{\top} \right\}\right]-\delta_{i}\frac{1}{n} \mathbf{{Z}}_{1} \mathbf{Z}_{2}^{\top}\mathbf{M}^{-1}(\hat\lambda_{i})\mathbf{M}^{-1}(\theta_{i})\left(\frac{1}{n}\mathbf{Z}_{2}\mathbf{Z}_{2}^{\top}\right)^{-1}\frac{1}{T}\mathbf{X}_{2}\mathbf{{X}}_{1}^{\top} \nonumber\\
&\quad+\theta_{i}\delta_{i} \frac{1}{n} \mathbf{{Z}}_{1} \mathbf{C}(\theta_{i})\mathbf{{X}}_{1}^{\top}  \nonumber\\
&=:\theta_{i}\Xi_{C1}-{\delta_{i}}\Xi_{C2}+\theta_{i}\delta_{i}\Xi_{C3}.
\end{align}

Similarly, we can write the last term $\Xi_{D}$ as
\begin{align}\label{seq41}
\Xi_{D}&=\theta_{i}  \frac{1}{T}\left[ \mathbf{{X}}_{1} \mathbf{D}(\theta_{i})\mathbf{{Z}}_{1}^{\top}-{\rm E}\left\{\mathbf{{X}}_{1} \mathbf{D}(\theta_{i})\mathbf{{Z}}_{1}^{\top} \right\}\right] \nonumber -\delta_{i} \frac{1}{T} \mathbf{{X}}_{1}
\mathbf{X}_{2}^{\top}\mathbf{M}^{-1}(\hat\lambda_{i})\mathbf{M}^{-1}(\theta_{i})\left(\frac{1}{n}\mathbf{Z}_{2}\mathbf{Z}_{2}^{\top}\right)^{-1}\frac{1}{n}\mathbf{Z}_{2}\mathbf{{Z}}_{1}^{\top} \nonumber\\
&\quad +\theta_{i} \delta_{i}\frac{1}{T} \mathbf{{X}}_{1} \mathbf{D}(\theta_{i})\mathbf{{Z}}_{1}^{\top}\nonumber\\
&=:\theta_{i}\Xi_{D1}-{\delta_{i}}\Xi_{D2}+\theta_{i}\delta_{i}\Xi_{D3}.
\end{align}
Putting (\ref{seq111})-(\ref{seq41}) into (\ref{object909909}), we have
 \begin{align}\label{object1}
\det(\theta_{i}\Theta_{1n}+\theta_{i}\delta_{i}\Theta_{2n}+\theta_{i}\Theta_{3n})=0,
\end{align}
where
\begin{align}
\Theta_{1n}&:=  (1+\delta_{i})^{2}\Xi_{A1}-\theta_{i}^{-1}\Xi_{B1}+\Xi_{C1}+\Xi_{D1},\label{nota1}\\
\Theta_{2n}&:= (1+\delta_{i})\Xi_{A2}
-(1+\delta_{i})\Xi_{A3}+\frac{1}{\theta_{i}\hat\lambda_{i}}\Xi_{B2}-\theta_{i}^{-1}\Xi_{C2}+\Xi_{C3}-\theta_{i}^{-1}\Xi_{D2}+\Xi_{D3},\label{nota2} \\
\Theta_{3n}&:=(1+\delta_{i})^{2}\Xi_{A4}-\theta_{i}^{-1}\Xi_{B3}.\label{nota3}
\end{align}
Multiplying both sides of the matrix in (\ref{object1}) by $\theta_i^{-1/2}\mathbf{U}$ from the left hand side and $\theta_i^{-1/2}\mathbf{U}^{\top}$ from the right hand side, we get
\begin{align}\label{object2}
\det\left\{\mathbf{U} ( \Theta_{1n}+\delta_{i}\Theta_{2n}+\Theta_{3n}) \mathbf{U}^{\top}\right\}=0.
\end{align}

Recall that $\mathbf{e}_{i}$ is the $q$-dimensional vector whose $i$-th element is $1$ and others are $0$.
 By Lemma \ref{lemma1.1} below, we have
\begin{align}\label{object43}
\sqrt{p}\hat{S}_{i}:=\sqrt{p}\mathbf{e}_{i}^{\top} \mathbf{U} \Theta_{1n}\mathbf{U}^{\top}\mathbf{e}_{i} \xrightarrow{d} \mathcal N(0,\tilde\sigma_{i}^{2}),
\end{align}
where $\tilde\sigma_{i}^{2}=\left(y+c\right)\left(1-y\right)^2\nu_i-y\left(1-y\right)\left(1-3y\right)-c\left(1-y\right)^2$.
It follows by Lemma \ref{lemma1.10} below  that
\begin{align}\label{object4600}
\|\mathbf{U}\Theta_{1n}\mathbf{U}^{\top}\|_{\infty}={\rm O}_{p}\left( \frac{q}{\sqrt{n}}+ \frac{\sum_j\lambda_j}{\sqrt{n}\lambda_i}\right).
\end{align}
By Lemma \ref{fisherlemma3} below, we also have
\begin{gather}
\max_{1\leq j\leq q}\left|\textbf{e}_j^{\top}\textbf{U} \Theta_{2n}\textbf{U}^{\top} \textbf{e}_j-(y-1)\right|={\rm O}_p\left(\frac{\sqrt q \delta_i}{\lambda_i}+\frac{\sqrt q}{\sqrt n}+\frac{\sqrt{\sum_j\lambda_j^2}}{\lambda_i^2}+\frac{\sqrt{\sum_j\lambda_j}}{\lambda_i}\right),\label{term2a}\\
\max_{1\leq j_1\neq j_2\leq q}\left|\textbf{e}_{j_1}^{\top}\textbf{U} \Theta_{2n}\textbf{U}^{\top} \textbf{e}_{j_2}\right|={\rm O}_p\left(\frac{ q\delta_i}{\lambda_i}+\frac{q}{\sqrt n}+\frac{{\sum_j\lambda_j}}{\lambda_i^2}+\frac{\sqrt{q\sum_j\lambda_j}}{\lambda_i}\right).\label{term2b}
\end{gather}
For the term $ \mathbf{U}\Theta_{3n}\mathbf{U}^{\top}$ in \eqref{object2}, by considering its $(j_1,j_2)$ entry for all $1\le j_1,j_2\le q$, we can easily get that
\begin{gather}
(1+\delta_{i})^{2}\mathbf{U}\Xi_{A4}\mathbf{U}^{\top}
=(1+\delta_{i})^{2}\left[1-\frac{p-q}{n}{\rm E}\left\{\tilde{m}_{\theta_{i}}(1) \right\}\right]\mathbf{I}_q\label{cc1},\\
\mathbf{U}\Xi_{B3}\mathbf{U}^{\top}=\left(1+\frac{p-q}{T}\left[-1+{\rm E}\left\{\tilde{ m}_{\theta_{i}}(1)\right\}\right]\right)\mathbf{\Lambda}_1.\label{cc2}
\end{gather}
By the definition of $\theta_{i}$ in (\ref{def_theta}),  we know
\begin{gather*}
1- \frac{p-q}{n}{\rm E}\left\{\tilde{ m}_{\theta_{i}}(1)\right\}=\frac{\lambda_{i} }{\theta_{i} }\left(1+\frac{p-q}{T}\left[-1+{\rm E}\left\{\tilde{ m}_{\theta_{i}}(1)\right\}\right]\right),
\end{gather*}
which,  together with the results in Lemma~\ref{fisherlemma2} below and Theorem~\ref{th1.2},
yields that
\begin{align}\label{cc3}
&(1+\delta_{i})^{2}\left[1- \frac{p-q}{n}{\rm E}\left\{\tilde{ m}_{\theta_{i}}(1)\right\} \right]-\frac{\lambda_{i} }{\theta_{i} }\left(1+\frac{p-q}{T}\left[-1+{\rm E}\left\{\tilde{ m}_{\theta_{i}}(1)\right\} \right]\right)\nonumber\\
=&2\delta_{i}\left[1- \frac{p-q}{n}{\rm E}\left\{\tilde{ m}_{\theta_{i}}(1)\right\} \right]+\delta_{i}^{2}\left[1- \frac{p-q}{n}{\rm E}\left\{\tilde{ m}_{\theta_{i}}(1)\right\} \right]= 2\delta_{i}\left\{1-\frac{p-q}{n}+{\rm o}(1)\right\}.
\end{align}

Combining (\ref{cc1})-(\ref{cc3}) and the definition of $\Theta_{3n}$ in (\ref{nota3}), we can get that, for $1\leq j\leq q$,
\begin{align*}
\mathbf{e}_{j}^{\top}\mathbf{U}\Theta_{3n}\mathbf{U}^{\top}\mathbf{e}_{j}^{\top}
=\left\{(1+\delta_{i})^2-\frac{\lambda_{j}}{\lambda_{i}}\right\}\left[1- \frac{p-q}{n}{\rm E}\left\{\tilde{ m}_{\theta_{i}}(1)\right\} \right],
\end{align*}
which converges to zero if and only if $\lambda_{j}=\lambda_{i}$ because  $(1+\delta_{i})^2-\lambda_{j}/\lambda_{i}>C>0$ for some constant $C$ if $\lambda_j\neq \lambda_i$ under Assumption~\ref{A4}. When $\lambda_{j}=\lambda_{i}$, we have
\begin{align}\label{object421}
\mathbf{e}_{j}^{\top} \mathbf{U}\Theta_{3n}\mathbf{U}^{\top}\mathbf{e}_{j}^{\top}=2\delta_{i}\left\{1-\frac{p-q}{n}+{\rm o}(1)\right\}.
\end{align}
Note that all off-diagonal entries of the matrix $ \mathbf{U}\Theta_{3n}\mathbf{U}^{\top}$ is zero, i.e.
\begin{align}\label{object422}
\mathbf{e}_{j_{1}}^{\top}\mathbf{U}\Theta_{3n}\mathbf{U}^{\top}\mathbf{e}_{j_{2}}=0, \forall 1\leq j_{1}\neq j_{2}\leq q.
\end{align}

Inserting (\ref{object43}), (\ref{object4600}), (\ref{term2a}), (\ref{term2b}), (\ref{object421}) and (\ref{object422}) into (\ref{object2}),
we can solve the determinant equation (\ref{object2}) and get the limiting distribution of $\delta_{i}(1\leq i \leq q)$ immediately.
Since diagonal elements of $\mathbf{U}\Theta_{3n}\mathbf{U}^{\top}$ are at least constant order, when $\mathbf{e}_{j}^{\top}\mathbf{U}\Theta_{3n}\mathbf{U}^{\top}\mathbf{e}_{j}^{\top}$ goes to infinity for some $j$'s, we can divide these rows by $\mathbf{e}_{j}^{\top}\mathbf{U}\Theta_{3n}\mathbf{U}^{\top}\mathbf{e}_{j}^{\top}$. In this way, we can get
\begin{align*}
\det \begin{pmatrix}{\rm O}_{p}(1)&\ldots &{\rm O}_{p}\left(*\right)&\ldots &{\rm O}_{p}\left(*\right)\\
\vdots&\ddots &\vdots&\ddots &\vdots\\
{\rm O}_{p}\left(*\right)&\ldots &\hat{S}_{i}+(1-y+{\rm o}_{p}(1))\delta_{i}&\ldots &{\rm O}_{p}\left(*\right)\\
\vdots&\ddots &\vdots&\ddots &\vdots\\
{\rm O}_{p}\left(*\right)&\ldots &{\rm O}_{p}\left(*\right)&\ldots &{\rm O}_{p}(1)
\end{pmatrix}=0
\end{align*}
where $\sqrt{p}\hat{S}_{i}\xrightarrow{d} \mathcal N(0,\tilde\sigma_{i}^{2})$ and
\begin{align*}
*=\frac{q}{\sqrt n}+\frac{{\sum_j\lambda_j}}{\sqrt n\lambda_i}+\frac{ q\delta_i^2}{\lambda_i}+\frac{\delta_i\sum_j\lambda_j}{\lambda_i^2}+\frac{\delta_i\sqrt{q\sum_j\lambda_j}}{\lambda_i}.
\end{align*}
By Leibniz formula for determinants, we can get that
$\hat{S}_{i}+ \left\{1-y+{\rm o}_{p}(1)\right\}\delta_{i}+q{\rm O}_{p}\left(*^2\right)=0$, that is
\begin{align*}
\hat{S}_{i}+ \left\{1-y+{\rm o}_{p}(1)\right\}\delta_{i}+{\rm O}_{p}\left(\frac{q^3}{ n}+\frac{{q(\sum_j\lambda_j})^2}{n\lambda_i^2}+\frac{ q^3\delta_i^4}{\lambda_i^2}+\frac{q\delta_i^2(\sum_j\lambda_j)^2}{\lambda_i^4}+\frac{q^2\delta_i^2{\sum_j\lambda_j}}{\lambda_i^2}\right)=0.
\end{align*}
Under Assumptions~\ref{A1} and \ref{A2}(a), we have $q={\rm o}(n^{\frac{1}{6}})$ and $\lambda_i^{-1}\sum_j\lambda_j={\rm o}(q^{-\frac{1}{2}}n^{\frac{1}{4}})$, then it follows that
\begin{gather*}
\frac{q^3}{ n}={\rm o}(n^{-\frac{1}{2}}),\quad \frac{{q^2\sum_j\lambda_j}}{n\lambda_i^2}={\rm o}(n^{-\frac{1}{2}}),\quad \frac{ q^3\delta_i^4}{\lambda_i^2}={\rm o}_p(\delta_i^2 n^{\frac{1}{2}}),\\
\frac{q\delta_i^2(\sum_j\lambda_j)^2}{\lambda_i^4}={\rm o}_p(\delta_i^2 n^{\frac{1}{2}}), \quad  \frac{q^2\delta_i^2{\sum_j\lambda_j}}{\lambda_i^2}={\rm o}_p(\delta_i^2 n^{\frac{1}{2}}).
\end{gather*}
It leads to
\begin{align*}
\hat{S}_{i}+ \left\{1-y+{\rm o}_{p}(1)\right\}\delta_{i}+{\rm o}_p(\delta_i^2 n^{\frac{1}{2}})+{\rm o}(n^{-\frac{1}{2}})=0.
\end{align*}
By multiplying $\sqrt p$ on both sides, we further obtain that
\begin{align*}
\sqrt{p}\hat{S}_{i}+ \left\{1-y+{\rm o}_{p}(1)\right\}\cdot\sqrt{p}\delta_{i}+{\rm o}_p\left(1\right)\cdot p\delta_i^2+{\rm o}\left(1\right)=0.
\end{align*}
Recalling that $\sqrt{p}\hat{S}_{i}\xrightarrow{d} \mathcal N(0,\tilde\sigma_{i}^{2})$, we can reach to
$\sqrt{p}\delta_{i} \xrightarrow{d} \mathcal N(0,\sigma_i^{2})$, where
\begin{align*}
\sigma_i^{2}=\frac{\tilde\sigma_{i}^{2}}{ (1-y)^{2}}=(y+c)\nu_i-c-\frac{y(1-3y)}{1-y}.
\end{align*}


Instead, we can consider the determinant
\begin{align}\label{caseb1}
\det\left\{ \tilde{\mathbf{\Lambda}}^{-\frac{1}{2}}\mathbf{U}( \theta_{i}\Theta_{1n}+\theta_{i}\delta_{i}\Theta_{2n}+\theta_{i}\Theta_{3n})\mathbf{U}^{\top}\tilde{\mathbf{\Lambda}}^{-\frac{1}{2}} \right\}=0,
\end{align}
where $\tilde{\mathbf{\Lambda}}=\text{diag}( \theta_{1},\ldots,\theta_{q})\in \mathbb{R}^{ q\times q}.$ Repeating all the derivations above, we can get
\begin{gather}
\|\tilde{\mathbf{\Lambda}}^{-\frac{1}{2}}\mathbf{U}\theta_i \Theta_{1n}\mathbf{U}^{\top}\tilde{\mathbf{\Lambda}}^{-\frac{1}{2}}\|_{\infty}={\rm O}_{p}\left( \frac{q}{\sqrt{n}}+ \frac{\lambda_i\sum_j\lambda_j^{-1}}{\sqrt{n}}\right),\label{caseb2}\\
\max_{1\leq j\leq q}\left|\textbf{e}_j^{\top}\tilde{\mathbf{\Lambda}}^{-\frac{1}{2}}\theta_{i}\mathbf{U} \Theta_{2n}\mathbf{U}^{\top} \tilde{\mathbf{\Lambda}}^{-\frac{1}{2}}\textbf{e}_j-(y-1)\frac{\theta_{i}}{\theta_{j}}\right|={\rm O}_p\left(\delta_i\sqrt{\sum_j\lambda_j^{-2}}+\frac{\lambda_i\sqrt{\sum_j\lambda_j^{-2}}}{\sqrt n}+\frac{\sqrt q}{\lambda_i}+\sqrt{\sum_j\lambda_j^{-1}}\right),\label{caseb3}\\
\max_{1\leq j_1\neq j_2\leq q}\left|\textbf{e}_{j_1}^{\top}\tilde{\mathbf{\Lambda}}^{-\frac{1}{2}}\theta_{i}\mathbf{U} \Theta_{2n}\mathbf{U}^{\top} \tilde{\mathbf{\Lambda}}^{-\frac{1}{2}}\textbf{e}_{j_2}\right|={\rm O}_p\left(\delta_i\sum_j\lambda_j^{-1}+\frac{\lambda_i\sum_j\lambda_j^{-1}}{\sqrt n}+\frac{q}{\lambda_i}+\sqrt{q\sum_j\lambda_j^{-1}}\right),\label{caseb4}\\
\mathbf{e}_{j}^{\top}\tilde{\mathbf{\Lambda}}^{-\frac{1}{2}}\theta_{i}\mathbf{U}\Theta_{3n}\mathbf{U}^{\top}\tilde{\mathbf{\Lambda}}^{-\frac{1}{2}}\mathbf{e}_{j}^{\top}
=\left\{(1+\delta_{i})^2-\frac{\lambda_{j}}{\lambda_{i}}\right\}\left[1- \frac{p-q}{n}{\rm E}\left\{\tilde{ m}_{\theta_{i}}(1)\right\} \right],\label{caseb5}\\
\mathbf{e}_{j_{1}}^{\top}\tilde{\mathbf{\Lambda}}^{-\frac{1}{2}}\theta_{i}\mathbf{U}\Theta_{3n}\mathbf{U}^{\top}\tilde{\mathbf{\Lambda}}^{-\frac{1}{2}}\mathbf{e}_{j_{2}}=0,\ \forall 1\leq j_{1}\neq j_{2}\leq q.\label{caseb6}
\end{gather}
Inserting \eqref{caseb2}-\eqref{caseb6} into \eqref{caseb1}, we can similarly prove $\sqrt{p}\delta_{i} \xrightarrow{d} \mathcal N(0,\sigma_i^{2})$ under Assumption~\ref{A2}(b).
Thus the proof is completed.\quad $\square$

\vspace{0.2in}

\noindent{\bf {Proof of Theorem \ref{th4.1}.}}
The proof of Theorem \ref{th4.1} is similar to that of Theorem \ref{th1.3}, the only difference is that we  take the $J_{i}\times J_{i}$ block as a typical object to analyse,
some useful lemmas can also be obtained from Lemmas \ref{lemma1.1}-\ref{fisherlemma3} below.
Similar arguments for deriving the proof of Theorem~4.1 in \cite{wang2017extreme} can be used. Thus, we omit the details.\quad $\square$


\section{Some Technical Lemmas}\label{sec4}
\begin{lemma}\label{fisherlemma2}
Suppose that Assumptions \ref{A1} and \ref{A3}  hold. For any
$\theta\rightarrow\infty$, we have $\widetilde{ m}_{\theta}(1)-1={\rm O}_{a.s.}({\theta}^{-1})$.
\end{lemma}
\noindent{\bf {Proof of lemma \ref{fisherlemma2}.}}
By the definition of $\widetilde{m}_{\theta}(z)$ in (\ref{def_m}),
\begin{align*}
\widetilde{ m}_{\theta}(1)&=\frac{1}{p-q}  \text{tr} \left(\textbf{I}_{p-q}-\frac{\textbf{F}_{0}}{\theta}  \right)^{-1}
=1+\frac{1}{p-q}  \text{tr} \left\{ \frac{\textbf{F}_{0}}{\theta}  \left(\textbf{I}_{p-q}-\frac{\textbf{F}_{0}}{\theta}  \right)^{-1}\right\},
\end{align*}
 we have
\begin{align*}
\widetilde{ m}_{\theta}(1)-1&=\frac{1}{p-q}\text{tr}\left\{\frac{\textbf{F}_{0}}{\theta}  \left(\textbf{I}_{p-q}-\frac{\textbf{F}_{0}}{\theta}  \right)^{-1}\right\} =\theta^{-1}\left(\frac{1}{p-q} \sum_{1\le j\le p-q}\frac{\mu_j}{1-\mu_j/\theta}\right).
\end{align*}
 Since all the eigenvalues of $\mathbf{F}_{0}$, namely $\mu_1\ge \ldots\ge \mu_{p-q}$, are almost surely bounded, we can get that $\widetilde{ m}_{\theta}(1)-1={\rm O}_{a.s.}({\theta}^{-1})$.\quad $\square$\\

Recall that $\mathbf{e}_{i}$ is the $q$-dimensional vector whose $i$-th element is $1$ and others are $0$, $\mathbf{U}^{\top}=(\mathbf{u}_1,\mathbf{u}_2,\ldots,\mathbf{u}_q)$, where $\mathbf{u}_i\in \mathbb{R}^q$ is the $i-th$ column of the matrix $\mathbf{U}^{\top}$. Then we  get the following lemma.

\begin{lemma}\label{lemma1.1}
For any fixed $1\leq i \leq q$, denote $\mathbf{G}_{ni}=\sqrt{p} \mathbf{U}\Theta_{1n}\mathbf{U}^{\top}$.
Under the assumptions of Theorem \ref{th1.3}, we have
\begin{align*}
\mathbf{e}^{\top}_{i}\mathbf{G}_{ni}\mathbf{e}_{i} \xrightarrow{d} \mathcal N(0, \tilde\sigma_{i}^{2}),
\end{align*}
where $\tilde\sigma^2_i=\left(y+c\right)\left(1-y\right)^2\nu_i-y\left(1-y\right)\left(1-3y\right)-c\left(1-y\right)^2$ and $\nu_{i}=\mathbb{E}|\mathbf{u}_i^{\top}\mathbf{Z}_1\mathbf{e}_1|^{4}$ for $1\leq i\leq q$.
\end{lemma}

\noindent{\bf {Proof of Lemma \ref{lemma1.1}.}} 
From the definition of $\Theta_{1n}$ in (\ref{nota1})
and the fact that $\mathbf{Y}_{1}=\Sigma_{1}^{-\frac{1}{2}}\mathbf{X}_{1}=\textbf{U}^{\top} {\mathbf{\Lambda}}^{-\frac{1}{2}}\textbf{U}\mathbf{X}_{1}$,
we have the decomposition
 \begin{align}\label{white0}
\mathbf{e}^{\top}_{i}\mathbf{G}_{ni}\mathbf{e}_{i}
=& \mathbf{u}_i^{\top}\Big[ \frac{(1+\delta_{i})^{2}\sqrt{p}}{n} \mathbf{Z}_{1}\left\{\textbf{I}_{n}-\mathbf{A}(\theta_{i})\right\}\mathbf{Z}_{1}^{\top}-
\frac{\lambda_{i}}{\theta_{i}}\frac{\sqrt{p}}{T} \mathbf{Y}_{1} \left\{\textbf{I}_{T}+\mathbf{B}(\theta_{i})\right\}\mathbf{Y}_{1}^{\top}\nonumber\\
&+\frac{\sqrt{p \lambda_{i}}}{n} \mathbf{Z}_{1} \mathbf{C}(\theta_{i})\mathbf{Y}_{1}^{\top}+ \frac{\sqrt{p \lambda_{i}}}{T}\mathbf{Y}_{1} \mathbf{D}(\theta_{i})\mathbf{Z}_{1}^{\top}\Big]\mathbf{u}_i
-{\rm E}(\cdot),
\end{align}
where ${\rm E}[\cdot]$ is the expectation of all the preceding terms after the  equal sign.

By Theorem \ref{th1.2}, $\delta_{i}$ converges in probability to $0$, thus we only need to consider the limit of
\begin{align*}
\mathbf{e}^{\top}_{i}\mathbf{\widetilde G}_{ni}\mathbf{e}_{i}
:= &\mathbf{u}_i^{\top}\Big[ \frac{\sqrt{p}}{n} \mathbf{Z}_{1}\left\{\textbf{I}_{n}-\mathbf{A}(\theta_{i})\right\}\mathbf{Z}_{1}^{\top}-
\frac{\lambda_{i}}{\theta_{i}}\frac{\sqrt{p}}{T }  \mathbf{Y}_{1} \left\{\textbf{I}_{T}+\mathbf{B}(\theta_{i})\right\}\mathbf{Y}_{1}^{\top}\nonumber+\frac{\sqrt{p \lambda_{i}}}{n} \mathbf{Z}_{1} \mathbf{C}(\theta_{i})\mathbf{Y}_{1}^{\top}\\
&+ \frac{\sqrt{p \lambda_{i}}}{T}\mathbf{Y}_{1} \mathbf{D}(\theta_{i})\mathbf{Z}_{1}^{\top}\Big]\mathbf{u}_i
-{\rm E}[\cdot].
\end{align*}

For the first two terms,  Theorem~7.2 in \cite{bai2008central} implies that,  for any $1\leq i\leq q$,
\begin{gather*}
\frac{1}{\sqrt n} \left[\mathbf{u}_i^{\top}\mathbf{Z}_{1}\left\{\textbf{I}_{n}-\mathbf{A}(\theta_{i})\right\}\mathbf{Z}_{1}^{\top}\mathbf{u}_i-\text{tr}\left\{\textbf{I}_{n}-\mathbf{A}(\theta_{i})\right\}\right]\xrightarrow{d} \mathcal N(0, \tilde\sigma_{i\textbf{A}}^{2}),\\
\frac{1}{\sqrt T} \left[\mathbf{u}_i^{\top}\mathbf{Y}_{1}\left\{\textbf{I}_{T}+\mathbf{B}(\theta_{i})\right\}\mathbf{Y}_{1}^{\top}\mathbf{u}_i-\text{tr}\left\{\textbf{I}_{T}+\mathbf{B}(\theta_{i})\right\}\right]\xrightarrow{d} \mathcal N(0, \tilde\sigma_{i\textbf{B}}^{2}),
\end{gather*}
with $\tilde\sigma_{iA}^{2}=\omega_{\textbf{I}_n-\textbf{A}(\theta_{i})}(\nu_{i}-3)+2\beta_{\textbf{I}_n-\textbf{A}(\theta_{i})}$ and $\tilde\sigma_{iB}^{2}=\omega_{\textbf{I}_T+\textbf{B}(\theta_{i})}(\nu_{i}-3)+2\beta_{\textbf{I}_T+\textbf{B}(\theta_{i})}$,
where
\begin{align*}
\nu_{i}&={\rm E}|\mathbf{u}_i^{\top}\mathbf{Z}_1\mathbf{e}_1|^{4}={\rm E}|\mathbf{u}_i^{\top}\mathbf{Y}_1\mathbf{e}_1|^{4},\\
\omega_{\textbf{I}_n-\textbf{A}(\theta_{i})}&=\lim_{n\rightarrow\infty}\frac{1}{n}\sum_{1\leq k\leq n}\left[\left\{\textbf{I}_n-\textbf{A}(\theta_{i})\right\}(k,k)\right]^2,\\
\beta_{\textbf{I}_n-\textbf{A}(\theta_{i})}&=\lim_{n\rightarrow\infty}\frac{1}{n}\text{tr}\left\{\textbf{I}_n-\textbf{A}(\theta_{i})\right\}^2,
\end{align*}
$\omega_{\textbf{I}_T+\textbf{B}(\theta_{i})}$ and $\beta_{\textbf{I}_T+\textbf{B}(\theta_{i})}$ are similarly defined. Here the fact that ${\rm E}|\mathbf{u}_i^{\top}\mathbf{Z}_1\mathbf{e}_1|^{4}={\rm E}|\mathbf{u}_i^{\top}\mathbf{Y}_1\mathbf{e}_1|^{4}$ is implied by Assumption~\ref{A3}.
Based on the facts that
\begin{align*}
&{\rm E}\left[\mathbf{u}_i^{\top}\mathbf{Z}_{1}\left\{\textbf{I}_{n}-\mathbf{A}(\theta_{i})\right\}\mathbf{Z}_{1}^{\top}\mathbf{u}_i\right]
={\rm E}\left(\text{tr}\left[\mathbf{Z}_{1}^{\top}\mathbf{u}_i\mathbf{u}_i^{\top}\mathbf{Z}_{1}\left\{\textbf{I}_{n}-\mathbf{A}(\theta_{i})\right\}\right]\right)\\
=&\text{tr}\left[{\rm E}\left(\mathbf{Z}_{1}^{\top}\mathbf{u}_i\mathbf{u}_i^{\top}\mathbf{Z}_{1}\right){\rm E}\left\{\textbf{I}_{n}-\mathbf{A}(\theta_{i})\right\}\right]
={\rm E}\left[\text{tr}\left\{\textbf{I}_{n}-\mathbf{A}(\theta_{i})\right\}\right]
=n-(p-q){\rm E}\left\{\widetilde m_{\theta_{i}}(1)\right\},
\end{align*}
and that $\widetilde m_{\theta_{i}}(1)-{\rm E}\left\{\widetilde m_{\theta_{i}}(1)\right\}={\rm O}_p(n^{-1})$,
we can get that
\begin{align*}
\frac{1}{\sqrt n}\left({\rm E}\left[\mathbf{u}_i^{\top}\mathbf{Z}_{1}\left\{\textbf{I}_{n}-\mathbf{A}(\theta_{i})\right\}\mathbf{Z}_{1}^{\top}\mathbf{u}_i^{\top}\right]-\text{tr}\left\{\textbf{I}_{n}-\mathbf{A}(\theta_{i})\right\}\right)={\rm o}_p(1).
\end{align*}
Then it follows that
\begin{align*}
\frac{1}{\sqrt n} \left(\mathbf{u}_i^{\top}\mathbf{Z}_{1}\left\{\textbf{I}_{n}-\mathbf{A}(\theta_{i})\right\}\mathbf{Z}_{1}^{\top}\mathbf{u}_i
-{\rm E}\left[\mathbf{u}_i^{\top}\mathbf{Z}_{1}\left\{\textbf{I}_{n}-\mathbf{A}(\theta_{i})\right\}\mathbf{Z}_{1}^{\top}\mathbf{u}_i^{\top}\right]\right)\xrightarrow{d} \mathcal N(0, \tilde\sigma_{i\textbf{A}}^{2}),
\end{align*}
and similarly,
\begin{align*}
\frac{1}{\sqrt T} \left(\mathbf{u}_i^{\top}\mathbf{Y}_{1}\left\{\textbf{I}_{T}+\mathbf{B}(\theta_{i})\right\}\mathbf{Y}_{1}^{\top}\mathbf{u}_i
-{\rm E}\left[\mathbf{u}_i^{\top}\mathbf{Y}_{1}\left\{\textbf{I}_{T}+\mathbf{B}(\theta_{i})\right\}\mathbf{Y}_{1}^{\top}\mathbf{u}_i\right]\right)\xrightarrow{d} \mathcal N(0, \tilde\sigma_{i\textbf{B}}^{2}).
\end{align*}
For other two terms, by the same approach in the proof of Theorem~\ref{th1.2}, we have that
\begin{align*}
\mathbf{u}_i^{\top}\left\{ \frac{\sqrt{p \lambda_{i}}}{n} \mathbf{Z}_{1} \mathbf{C}(\theta_{i})\mathbf{Y}_{1}^{\top}+\frac{\sqrt{p \lambda_{i}}}{T}\mathbf{Y}_{1} \mathbf{D}(\theta_{i})\mathbf{Z}_{1}^{\top}\right\}\mathbf{u}_i={\rm O}_p\left(\frac{1}{\sqrt \lambda_i}\right).
\end{align*}
By all these arguments above, we can derive that $\mathbf{e}^{\top}_{i}\mathbf{G}_{ni}\mathbf{e}_{i}\xrightarrow{d} \mathcal N(0, \tilde\sigma_{i}^{2})$ with $\tilde\sigma_i^2=y\sigma_{i\textbf{A}}^2+c(1-y)^2\sigma_{i\textbf{B}}^2$.\\

We compute $\omega_{\textbf{I}_n-\textbf{A}(\theta_{i})}$, $\beta_{\textbf{I}_n-\textbf{A}(\theta_{i})}$, $\omega_{\textbf{I}_T+\textbf{B}(\theta_{i})}$ and $\beta_{\textbf{I}_T+\textbf{B}(\theta_{i})}$ in the following.
By the derivations in the proof of  Lemma 6 in \cite{wang2017extreme},
\begin{align*}
\left\{\textbf{I}_n-\textbf{A}(\theta_{i})\right\}(k,k)&=1-\left\{\mathbf{Z}_{2}^{\top}\mathbf{M}(\theta_{i})^{-1}\left(\frac{1}{n}\mathbf{Z}_{2}\mathbf{Z}_{2}^{\top}\right)^{-1}\frac{1}{n}\mathbf{Z}_{2}\right\} (k,k)\\
&=1-\frac{\theta_{i}}{n}\left\{\mathbf{Z}_{2}^{\top}\left(\theta_{i}\cdot\frac{1}{n}\mathbf{Z}_{2}\mathbf{Z}_{2}^{\top}
-\frac{1}{T}\mathbf{X}_{2}\mathbf{X}_{2}^{\top}\right)^{-1}\mathbf{Z}_{2}\right\} (k,k)\\
&=\frac{1}{1+\frac{\theta_{i}}{n}\left\{\eta_{k}^{\top}\left(\theta_{i}\frac{1}{n}\mathbf{Z}_{2k}\mathbf{Z}_{2k}^{\top}-\frac{1}{T}\mathbf{X}_{2}\mathbf{X}_{2}^{\top}  \right)^{-1}\eta_{k}\right\}},
\end{align*}
where $\eta_k$ is the $k$-th column of $\mathbf{Z}_2$ and $\mathbf{Z}_{2k}$ is defined by removing the $k$-th column of $\mathbf{Z}_2$.

Note that
\begin{align}
&\left(\frac{1}{n}\mathbf{Z}_{2k}\mathbf{Z}_{2k}^{\top}-\frac{1}{\theta_{i}T}\mathbf{X}_{2}\mathbf{X}_{2}^{\top}   \right)^{-1}-\left(\frac{1}{n}\mathbf{Z}_{2k}\mathbf{Z}_{2k}^{\top}\right)^{-1}\nonumber\\
=&\left(\frac{1}{n}\mathbf{Z}_{2i}\mathbf{Z}_{2i}^{\top}-\frac{1}{\theta_{i}T}\mathbf{X}_{2}\mathbf{X}_{2}^{\top}\right)^{-1}\left\{\frac{1}{n}\mathbf{Z}_{2k}\mathbf{Z}_{2k}^{\top}-\left(\frac{1}{n}\mathbf{Z}_{2k}\mathbf{Z}_{2k}^{\top}-\frac{1}{\theta_{i}T}\mathbf{X}_{2}\mathbf{X}_{2}^{\top}\right)\right\}\left(\frac{1}{n}\mathbf{Z}_{2k}\mathbf{Z}_{2k}^{\top}\right)^{-1}\nonumber\\
=&\theta_{i}^{-1}\left(\frac{1}{n}\mathbf{Z}_{2i}\mathbf{Z}_{2i}^{\top}-\frac{1}{\theta_{i}T}\mathbf{X}_{2}\mathbf{X}_{2}^{\top}\right)^{-1}\left(\frac{1}{T}\mathbf{X}_{2}\mathbf{X}_{2}^{\top}\right)\left(\frac{1}{n}\mathbf{Z}_{2k}\mathbf{Z}_{2k}^{\top}\right)^{-1}
\end{align}
and
\begin{align}
\frac{1}{p-q-1}{\rm tr}\left(\frac{1}{n}\mathbf{Z}_{2k}\mathbf{Z}_{2k}^{\top}\right)^{-1}= \mathcal{S}_{\rm MP}\left(0\right)+{\rm O}_p\left(p^{-1}\right)=\frac{1}{1-y}+{\rm O}_p\left(p^{-1}\right),
\end{align}
where $\mathcal{S}_{\rm MP}$ denotes the Stieltjes transform of the Marcenko-Pastur law.
Then we have that
\begin{align*}
\frac{1}{p-q}\theta_{i}{\rm E}\left\{\text{tr}\left(\theta_{i}\frac{1}{n}\mathbf{Z}_{2k}\mathbf{Z}_{2k}^{\top}-\frac{1}{T}\mathbf{X}_{2}\mathbf{X}_{2}^{\top} \right)^{-1}\right\}
&={\rm E}\left\{\frac{1}{p-q}\text{tr}\left(\frac{1}{n}\mathbf{Z}_{2k}\mathbf{Z}_{2k}^{\top}\right)^{-1}+{\rm O}_{a.s.}\left(\theta_i^{-1}\right)\right\}\nonumber\\
&={\rm E}\left\{\frac{1}{1-y}+{\rm O}_{a.s.}\left(\theta_i^{-1}\right)+{\rm O}_p\left(p^{-1}\right)\right\}\rightarrow \frac{1}{1-y}.
\end{align*}
By Lemma A.2. in \cite{wang2017extreme}, it holds that
\begin{align*}
\left\{\textbf{I}_{n}-\mathbf{A}(\theta_{i})\right\}(k,k)\rightarrow \frac{1}{1+y(1-y)^{-1}}=1-y,
\end{align*}
which implies
\begin{align*}
\omega_{\textbf{I}_n-\textbf{A}(\theta_{i})}=\lim_{n\rightarrow\infty}\frac{1}{n}\sum_{1\leq k\leq n}\left[\left\{\textbf{I}_n-\mathbf{A}(\theta_{i})\right\}(k,k)\right]^2=(1-y)^2.
\end{align*}
By the similar argument, we can obtain that
\begin{align*}
\omega_{\textbf{I}_T+\textbf{B}(\theta_{i})}=\lim_{T\rightarrow\infty}\frac{1}{T} \sum_{1\leq k\leq T} \left[\left\{\textbf{I}_{T}+\mathbf{B}(\theta_{i})\right\}(k,k)\right]^{2}=1.
\end{align*}
Now we come to the calculation of $\beta_{\textbf{I}_n-\textbf{A}(\theta_{i})}$ and $\beta_{\textbf{I}_T+\textbf{B}(\theta_{
i})}$. Since $\theta_i\rightarrow +\infty$ as $n$ goes to infinity, we have
\begin{align*}
&\lim_{n\rightarrow\infty}\int_{-\infty}^{\infty} \frac{\theta_{i}}{\theta_{i}-x}dF_n(x)=1,\  \lim_{n\rightarrow\infty}\int_{-\infty}^{\infty}\frac{\theta_{i}^2}{(\theta_{i}-x)^2}dF_n(x)=1,\\
&\lim_{T\rightarrow\infty}\int_{-\infty}^{\infty} \frac{x}{\theta_{i}-x}dF_n(x)=0,\
\lim_{T\rightarrow\infty}\int_{-\infty}^{\infty}\frac{x^2}{(\theta_{i}-x)^2}dF_n(x)=0.
\end{align*}
Then these calculations lead to
\begin{align*}
\beta_{\textbf{I}_n-\textbf{A}(\theta_{i})}&=\lim_{n\rightarrow\infty}\frac{1}{n}\text{tr}\left\{\textbf{I}_n-\textbf{A}(\theta_{i})\right\}^2=\lim_{n\rightarrow\infty}\frac{1}{n}\text{tr}\left\{\textbf{I}_n-2\textbf{A}(\theta_{i})+\textbf{A}^2(\theta_{i})\right\}\\
&=1-2\lim_{n\rightarrow\infty}\left(\frac{p-q}{n}\int_{-\infty}^{\infty} \frac{\theta_{i}}{\theta_{i}-x}dF_n(x)\right)+\lim_{n\rightarrow\infty}\left\{\frac{p-q}{n}\int_{-\infty}^{\infty}\frac{\theta_{i}^2}{(\theta_{i}-x)^2} dF_n(x)\right\}\\
&=1-2y+y=1-y,\\
\beta_{\textbf{I}_T+\textbf{B}(\theta_{i})}&=\lim_{T\rightarrow\infty}\frac{1}{T}\text{tr}\left\{\textbf{I}_T+\textbf{B}(\theta_{i})\right\}^2=\lim_{T\rightarrow\infty}\frac{1}{T}\text{tr}\left\{\textbf{I}_T+2\textbf{B}(\theta_{i})+\textbf{B}^2(\theta_{i})\right\}\\
&=1+2\lim_{T\rightarrow\infty}\left\{\frac{p-q}{T}\int_{-\infty}^{\infty} \frac{x}{\theta_{i}-x}dF_n(x)\right\}+\lim_{T\rightarrow\infty}\left\{\frac{p-q}{T}\int_{-\infty}^{\infty}\frac{x^2}{(\theta_{i}-x)^2} dF_n(x)\right\}\\
&=1+0+0=1.
\end{align*}
Thus, we can write
\begin{align*}
\tilde\sigma_i^2&=y\tilde\sigma_{i\textbf{A}}^2+c(1-y)^2\tilde\sigma_{i\textbf{B}}^2\\
&=y\{\omega_{\textbf{I}_n-\textbf{A}(\theta_{i})}(\nu_{i}-3)+2\beta_{\textbf{I}_n-\textbf{A}(\theta_{i})}\}+c(1-y)^2\{\omega_{\textbf{I}_T+\textbf{B}(\theta_{i})}(\nu_{i}-3)+2\beta_{\textbf{I}_T+\textbf{B}(\theta_{i})}\}\\
&=\left(y+c\right)\left(1-y\right)^2\nu_i-y\left(1-y\right)\left(1-3y\right)-c\left(1-y\right)^2.
\end{align*}
Thus the proof is completed.\quad$\square$

\vspace{0.2in}

\begin{lemma}\label{lemma1.10}
Under the assumptions of Theorem \ref{th1.3}, 
\begin{align*}
\| \mathbf{U}\Theta_{1n}\mathbf{U}^{\top} \|_{\infty}={\rm O}_{p}\left( \frac{q}{\sqrt{n}}+ \frac{\sum_j\lambda_j}{\sqrt{n}\lambda_i}\right).
\end{align*}
\end{lemma}

\noindent{\bf {Proof of Lemma \ref{lemma1.10}.}} By the definition of $\Theta_{1n}$ in (\ref{nota1}) again, we know
\begin{align}\label{new1}
\Theta_{1n}
=& (1+\delta_{i})^{2}\frac{1}{n} \mathbf{Z}_{1}\left\{\textbf{I}_{n}-\mathbf{A}(\theta_{i})\right\}\mathbf{Z}_{1}^{\top}-\theta_{i}^{-1}
\frac{1}{T}  \mathbf{X}_{1} \left\{\textbf{I}_{T}+\mathbf{B}(\theta_{i})\right\}\mathbf{X}_{1}^{\top}\nonumber\\
&+\frac{1}{n} \mathbf{Z}_{1} \mathbf{C}(\theta_{i})\mathbf{X}_{1}^{\top}+ \frac{1}{T}\mathbf{X}_{1} \mathbf{D}(\theta_{i})\mathbf{Z}_{1}^{\top}
-{\rm E}(\cdot),
\end{align}
where $\rm{E}(\cdot)$ is the expectation of all the preceding terms.

Denote
\begin{align*}
&\eta_{n1}=\frac{1}{n} \mathbf{Z}_{1}\left\{\textbf{I}_{n}-\mathbf{A}(\theta_{i})\right\}\mathbf{Z}_{1}^{\top}-{\rm E}\left[\frac{1}{n} \mathbf{Z}_{1}\left\{\textbf{I}_{n}-\mathbf{A}(\theta_{i})\right\}\mathbf{Z}_{1}^{\top}\right],\\
&\eta_{n2}=\frac{1}{T}\mathbf{Y}_{1} \left\{\textbf{I}_{T}+\mathbf{B}(\theta_{i})\right\}\mathbf{Y}_{1}^{\top}-{\rm E}\left[\frac{1}{T}\mathbf{Y}_{1} \left\{\textbf{I}_{T}+\mathbf{B}(\theta_{i})\right\}\mathbf{Y}_{1}^{\top}\right],\\
&\eta_{n3}=\sqrt{\theta_{i}}\frac{1}{n} \mathbf{Z}_{1} \mathbf{C}(\theta_{i})\mathbf{Y}_{1}^{\top}-{\rm E}\left\{\sqrt{\theta_{i}}\frac{1}{n} \mathbf{Z}_{1} \mathbf{C}(\theta_{i})\mathbf{Y}_{1}^{\top}\right\},\\
&\eta_{n4}=\sqrt{\theta_{i}}\frac{1}{T}\mathbf{Y}_{1} \mathbf{D}(\theta_{i})\mathbf{Z}_{1}^{\top}-{\rm E}\left\{\sqrt{\theta_{i}}\frac{1}{T}\mathbf{Y}_{1} \mathbf{D}(\theta_{i})\mathbf{Z}_{1}^{\top}\right\}.
\end{align*}
By the fact $\mathbf{X}_{1}=\mathbf{U}^{\top}\Lambda_1^{\frac{1}{2}}\mathbf{U}\mathbf{Y}_{1}$, we can write
\begin{align}\label{new2}
 \mathbf{U}\Theta_{1n}\mathbf{U}^{\top} :=\sum_{i=1}^{4}\textbf{V}_{ni},
\end{align}
where
\begin{align}
\textbf{V}_{n1}&=(1+\delta_{i})^{2} \mathbf{U}   \left(\frac{1}{n} \mathbf{Z}_{1}\left\{\textbf{I}_{n}-\mathbf{A}(\theta_{i})\right\}\mathbf{Z}_{1}^{\top}-{\rm E}\left[\frac{1}{n} \mathbf{Z}_{1}\left\{\textbf{I}_{n}-\mathbf{A}(\theta_{i})\right\}\mathbf{Z}_{1}^{\top}\right] \right)   \mathbf{U}^{\top}\nonumber \\
&= (1+\delta_{i})^{2} \mathbf{U}  \eta_{n1}   \mathbf{U}^{\top},\label{new3}\\
\textbf{V}_{n2}&=-\theta_{i}^{-1}\mathbf{U}   \left(\frac{1}{T}\mathbf{X}_{1} \left\{\textbf{I}_{T}+\mathbf{B}(\theta_{i})\right\}\mathbf{X}_{1}^{\top}-{\rm E}\left[\frac{1}{T}\mathbf{X}_{1} \left\{\textbf{I}_{T}+\mathbf{B}(\theta_{i})\right\}\mathbf{X}_{1}^{\top}\right]  \right)  \mathbf{U}^{\top} \nonumber\\
&=-\theta_{i}^{-1}\mathbf{\Lambda}^{\frac{1}{2}} \mathbf{U}   \left(\frac{1}{T}\mathbf{Y}_{1} \left\{\textbf{I}_{T}+\mathbf{B}(\theta_{i})\right\}\mathbf{Y}_{1}^{\top}-{\rm E}\left[\frac{1}{T}\mathbf{Y}_{1} \left\{\textbf{I}_{T}+\mathbf{B}(\theta_{i})\right\}\mathbf{Y}_{1}^{\top}\right]  \right)  \mathbf{U}^{\top}\mathbf{\Lambda}^{\frac{1}{2}}\nonumber \\
&=-\theta_{i}^{-1}\mathbf{\Lambda}^{\frac{1}{2}} \mathbf{U}   \eta_{n2}  \mathbf{U}^{\top}\mathbf{\Lambda}^{\frac{1}{2}},\label{new4}\\
\textbf{V}_{n3}&=\mathbf{U}\left[\frac{1}{n} \mathbf{Z}_{1} \mathbf{C}(\theta_{i})\mathbf{X}_{1}^{\top}-{\rm E}\left\{\frac{1}{n} \mathbf{Z}_{1} \mathbf{C}(\theta_{i})\mathbf{X}_{1}^{\top}\right\}\right]   \mathbf{U}^{\top}\nonumber\\
&=\mathbf{U}\left[\frac{1}{n} \mathbf{Z}_{1} \mathbf{C}(\theta_{i})\mathbf{Y}_{1}^{\top}-{\rm E}\left\{\frac{1}{n} \mathbf{Z}_{1} \mathbf{C}(\theta_{i})\mathbf{Y}_{1}^{\top}\right\}\right]    \mathbf{U}^{\top} \mathbf{\Lambda}^{\frac{1}{2}}\nonumber\\
&=\mathbf{U}\left[\frac{1}{n} \mathbf{Z}_{1} \mathbf{C}(\theta_{i})\mathbf{Y}_{1}^{\top}-{\rm E}\left\{\frac{1}{n} \mathbf{Z}_{1} \mathbf{C}(\theta_{i})\mathbf{Y}_{1}^{\top}\right\}\right] \mathbf{U}^{\top} \mathbf{\Lambda}^{\frac{1}{2}}\nonumber\\
&=\theta_{i}^{-\frac{1}{2}}\mathbf{U}\eta_{n3} \mathbf{U}^{\top}\mathbf{\Lambda}^{\frac{1}{2}}\label{new5}\\
\textbf{V}_{n4}&=\mathbf{U}\left[\frac{1}{T}\mathbf{X}_{1} \mathbf{D}(\theta_{i})\mathbf{Z}_{1}^{\top}-{\rm E}\left\{\frac{1}{T}\mathbf{X}_{1} \mathbf{D}(\theta_{i})\mathbf{Z}_{1}^{\top}\right\}\right]\mathbf{U}^{\top} =\theta_{i}^{-\frac{1}{2}}\mathbf{\Lambda}^{\frac{1}{2}}\mathbf{U}\eta_{n4} \mathbf{U}^{\top}.\label{new6}
\end{align}

Similarly as the arguments in  the proof of Lemma~\ref{lemma1.1}, it holds that, for $1\leq j_1,j_2\leq q$,
\begin{gather*}
\textbf{e}_{j_1}^{\top}\eta_{n1}\textbf{e}_{j_2}={\rm O}_p\left(\frac{1}{\sqrt n}\right),\quad \textbf{e}_{j_1}^{\top}\eta_{n2}\textbf{e}_{j_2}={\rm O}_p\left(\frac{1}{\sqrt n}\right),\\
\textbf{e}_{j_1}^{\top}\eta_{n3}\textbf{e}_{j_2}={\rm O}_p\left(\frac{1}{\sqrt {n\lambda_i}}\right),\quad
\textbf{e}_{j_1}^{\top}\eta_{n4}\textbf{e}_{j_2}={\rm O}_p\left(\frac{1}{\sqrt {n\lambda_i}}\right).
\end{gather*}

Noting  that $\mathbf{U}$ is an orthogonal matrix, we have that
\begin{align*}
&\mathbf{e}_{j_{1}}^{\top}\textbf{V}_{n1}\mathbf{e}_{j_{2}}=\mathbf{e}_{j_{1}}^{\top}(1+\delta_{i})^{2} \mathbf{U}   \eta_{n1}   \mathbf{U}^{\top} \mathbf{e}_{j_{2}}={\rm O}_{p}\left(\frac{1}{\sqrt{n}}\right),\\
&\mathbf{e}_{j_{1}}^{\top}\textbf{V}_{n2}\mathbf{e}_{j_{2}}=-\mathbf{e}_{j_{1}}^{\top}\theta_{i}^{-1}\mathbf{\Lambda}^{\frac{1}{2}} \mathbf{U}   \eta_{n2}  \mathbf{U}^{\top}\mathbf{\Lambda}^{\frac{1}{2}}\mathbf{e}_{j_{2}}=\lambda_i^{-1}\lambda_{j_1}^{\frac{1}{2}}\lambda_{j_2}^{\frac{1}{2}}\cdot {\rm O}_{p}\left(\frac{1}{\sqrt{n}}\right),\\
&\mathbf{e}_{j_{1}}^{\top}\textbf{V}_{n3}\mathbf{e}_{j_{2}}=\mathbf{e}_{j_{1}}^{\top}\theta_{i}^{-\frac{1}{2}}\mathbf{U}\eta_{n3} \mathbf{U}^{\top}\mathbf{\Lambda}^{\frac{1}{2}}\mathbf{e}_{j_{2}}=\lambda_i^{-1}\lambda_{j_2}^{\frac{1}{2}}\cdot {\rm O}_{p}\left(\frac{1}{\sqrt{n}}\right),\\
&\mathbf{e}_{j_{1}}^{\top}\textbf{V}_{n4}\mathbf{e}_{j_{2}}=\mathbf{e}_{j_{1}}^{\top}\theta_{i}^{-\frac{1}{2}}\mathbf{\Lambda}^{\frac{1}{2}}\mathbf{U}\eta_{n4} \mathbf{U}^{\top}\mathbf{e}_{j_{2}}=\lambda_i^{-1}\lambda_{j_1}^{\frac{1}{2}}\cdot {\rm O}_{p}\left(\frac{1}{\sqrt{n}}\right).
\end{align*}
Then by Chebyshev's inequality, we can deduce that
\begin{align*}
\| \textbf{V}_{n1} \|_{\infty}&={\rm O}_{p}\left( \frac{q}{\sqrt{n}}\right),\quad
\| \textbf{V}_{n2} \|_{\infty}={\rm O}_{p}\left( \frac{\sum_j\lambda_j}{\sqrt{n}\lambda_i}\right),\quad
\| \textbf{V}_{n3}+\textbf{V}_{n4} \|_{\infty}={\rm O}_{p}\left( \frac{\sqrt {q\sum_j\lambda_j}}{\sqrt{n}\lambda_i}\right),
\end{align*}
where $\sqrt {q\sum_j\lambda_j}={\rm o}(\sum_j\lambda_j) $.
Thus we complete the proof by (\ref{new2}).   \quad$\square$

\vspace{0.1in}

\begin{lemma}\label{fisherlemma3}
Under the assumptions of Theorem \ref{th1.3},
\begin{gather}
\max_{1\leq j\leq q}\left|\textbf{e}_j^{\top}\mathbf{U} \Theta_{2n}\mathbf{U}^{\top} \textbf{e}_j-(y-1)\right|={\rm O}_p\left(\frac{\sqrt q \delta_i}{\lambda_i}+\frac{\sqrt q}{\sqrt n}+\frac{\sqrt{\sum_j\lambda_j^2}}{\lambda_i^2}+\frac{\sqrt{\sum_j\lambda_j}}{\lambda_i}\right),\\
\max_{1\leq j_1\neq j_2\leq q}|\textbf{e}_{j_1}^{\top}\mathbf{U} \Theta_{2n}\mathbf{U}^{\top} \textbf{e}_{j_2}|={\rm O}_p\left(\frac{ q\delta_i}{\lambda_i}+\frac{q}{\sqrt n}+\frac{{\sum_j\lambda_j}}{\lambda_i^2}+\frac{\sqrt{q\sum_j\lambda_j}}{\lambda_i}\right).
\end{gather}
\end{lemma}
\noindent{\bf {Proof of Lemma~\ref{fisherlemma3}}.} Recall the definition of $\Theta_{2n}$ in \eqref{nota2}:
\begin{align*}
\Theta_{2n}=& (1+\delta_{i}) \frac{1}{n}\mathbf{Z}_{1}\mathbf{Z}_{2}^{\top}\mathbf{M}^{-1}(\hat\lambda_{i})\mathbf{M}^{-1}(\theta_{i})\left(\frac{1}{n}\mathbf{Z}_{2}\mathbf{Z}_{2}^{\top}\right)^{-1}\frac{1}{n}\mathbf{Z}_{2} \mathbf{Z}_{1}^{\top} -(1+\delta_{i}) \frac{1}{n} \mathbf{Z}_{1}\mathbf{Z}_{1}^{\top} \\
&+\frac{1}{\hat\lambda_{i} \theta_{i}}
\frac{1}{T} \mathbf{X}_{1}\mathbf{X}_{2}^{\top}\mathbf{M}^{-1}(\hat\lambda_{i})\mathbf{M}^{-1}(\theta_{i})\left(\frac{1}{n}\mathbf{Z}_{2}\mathbf{Z}_{2}^{\top}\right)^{-1}\frac{1}{T}\mathbf{X}_{2}
\mathbf{X}_{1}^{\top} \\
&-\frac{(1+\delta_{i})}{\hat\lambda_{i}} \frac{1}{n} \mathbf{Z}_{1}
\mathbf{Z}_{2}^{\top}\mathbf{M}^{-1}(\hat\lambda_{i})\mathbf{M}^{-1}(\theta_{i})\left(\frac{1}{n}\mathbf{Z}_{2}\mathbf{Z}_{2}^{\top}\right)^{-1}\frac{1}{T}\mathbf{X}_{2}\mathbf{X}_{1}^{\top} +\frac{1}{n} \mathbf{Z}_{1} \mathbf{C}(\theta_{i})\mathbf{X}_{1}^{\top}   \nonumber\\
&-\frac{(1+\delta_{i})}{\hat\lambda_{i}} \frac{1}{T} \mathbf{X}_{1}
\mathbf{X}_{2}^{\top}\mathbf{M}^{-1}(\hat\lambda_{i})\mathbf{M}^{-1}(\theta_{i})\left(\frac{1}{n}\mathbf{Z}_{2}\mathbf{Z}_{2}^{\top}\right)^{-1}\frac{1}{n}\mathbf{Z}_{2}\mathbf{Z}_{1}^{\top} + \frac{1}{T}\mathbf{X}_{1} \mathbf{D}(\theta_{i})\mathbf{Z}_{1}^{\top}.
\end{align*}
Noting that
\begin{align*}
\mathbf{M}^{-1}(\hat\lambda_{i})-\mathbf{M}^{-1}(\theta_{i})=\mathbf{M}^{-1}(\hat\lambda_{i})\left\{\mathbf{M}(\theta_{i})
-\mathbf{M}(\hat\lambda_{i})\right\}\mathbf{M}^{-1}(\theta_{i})=-\frac{\delta_{i}}{\hat\lambda_{i}}\mathbf{M}^{-1}(\hat\lambda_{i})\textbf{F}_0\mathbf{M}^{-1}(\theta_{i}),
\end{align*}
we decompose the first term in $\Theta_{2n}$ as
\begin{align*}
&\frac{1}{n}\mathbf{Z}_{1}\mathbf{Z}_{2}^{\top}\mathbf{M}^{-1}(\hat\lambda_{i})\mathbf{M}^{-1}(\theta_{i})\left(\frac{1}{n}\mathbf{Z}_{2}\mathbf{Z}_{2}^{\top}\right)^{-1}\frac{1}{n}\mathbf{Z}_{2} \mathbf{Z}_{1}^{\top}\\
=&\frac{1}{n}\mathbf{Z}_{1}\mathbf{Z}_{2}^{\top}\left\{\mathbf{M}^{-1}(\hat\lambda_{i})-\mathbf{M}^{-1}(\theta_{i})\right\}\mathbf{M}^{-1}(\theta_{i})\left(\frac{1}{n}\mathbf{Z}_{2}\mathbf{Z}_{2}^{\top}\right)^{-1}\frac{1}{n}\mathbf{Z}_{2} \mathbf{Z}_{1}^{\top}+\frac{1}{n}\mathbf{Z}_{1}\mathbf{Z}_{2}^{\top}\mathbf{M}^{-2}(\theta_{i})\left(\frac{1}{n}\mathbf{Z}_{2}\mathbf{Z}_{2}^{\top}\right)^{-1}\frac{1}{n}\mathbf{Z}_{2} \mathbf{Z}_{1}^{\top}\\
=&-\frac{\delta_{i}}{\hat\lambda_{i}}\frac{1}{n}\mathbf{Z}_{1}\mathbf{Z}_{2}^{\top}\mathbf{M}^{-1}(\hat\lambda_{i})\textbf{F}_0\mathbf{M}^{-2}(\theta_{i})\left(\frac{1}{n}\mathbf{Z}_{2}\mathbf{Z}_{2}^{\top}\right)^{-1}\frac{1}{n}\mathbf{Z}_{2} \mathbf{Z}_{1}^{\top}+\frac{1}{n}\mathbf{Z}_{1}\mathbf{Z}_{2}^{\top}\mathbf{M}^{-2}(\theta_{i})\left(\frac{1}{n}\mathbf{Z}_{2}\mathbf{Z}_{2}^{\top}\right)^{-1}\frac{1}{n}\mathbf{Z}_{2} \mathbf{Z}_{1}^{\top}.
\end{align*}
On one hand, similar to the arguments in the proof of Theorem~\ref{th1.2}, we can derive that
\begin{gather*}
\max_{1\leq j\leq q}\left|\textbf{e}_j^{\top}\frac{\delta_{i}}{\hat\lambda_{i}}\frac{1}{n}\mathbf{Z}_{1}\mathbf{Z}_{2}^{\top}\mathbf{M}^{-1}(\hat\lambda_{i})\textbf{F}_0\mathbf{M}^{-2}(\theta_{i})\left(\frac{1}{n}\mathbf{Z}_{2}\mathbf{Z}_{2}^{\top}\right)^{-1}\frac{1}{n}\mathbf{Z}_{2} \mathbf{Z}_{1}^{\top}\textbf{e}_j\right|={\rm O}_{p}\left(\frac{\sqrt q \delta_i}{\lambda_i}\right),\\
\max_{1\leq j_1\neq j_2\leq q}\left|\textbf{e}_{j_1}^{\top}\frac{\delta_{i}}{\hat\lambda_{i}}\frac{1}{n}\mathbf{Z}_{1}\mathbf{Z}_{2}^{\top}
\mathbf{M}^{-1}(\hat\lambda_{i})\textbf{F}_0\mathbf{M}^{-2}(\theta_{i})\left(\frac{1}{n}\mathbf{Z}_{2}\mathbf{Z}_{2}^{\top}\right)^{-1}\frac{1}{n}\mathbf{Z}_{2} \mathbf{Z}_{1}^{\top}\textbf{e}_{j_2}\right|
={\rm O}_p\left(\frac{ q\delta_i}{\lambda_i}\right).
\end{gather*}
On the other hand, similar to the proof of Lemma~\ref{lemma1.1}, we can get that
\begin{gather*}
\frac{1}{n}\left(\textbf{e}_j^{\top}\mathbf{Z}_{1}\mathbf{Z}_{2}^{\top}\mathbf{M}^{-2}(\theta_{i})\left(\frac{1}{n}\mathbf{Z}_{2}\mathbf{Z}_{2}^{\top}\right)^{-1}\frac{1}{n}\mathbf{Z}_{2} \mathbf{Z}_{1}^{\top}\textbf{e}_j-{\rm E}\left[\text{tr}\left\{\mathbf{M}^{-2}(\theta_{i})\right\}\right]\right)={\rm O}_p\left(\frac{1}{\sqrt n}\right),\\
\frac{1}{n}\textbf{e}_{j_1}^{\top}\mathbf{Z}_{1}\mathbf{Z}_{2}^{\top}\mathbf{M}^{-2}(\theta_{i})\left(\frac{1}{n}\mathbf{Z}_{2}\mathbf{Z}_{2}^{\top}\right)^{-1}\frac{1}{n}\mathbf{Z}_{2} \mathbf{Z}_{1}^{\top}\textbf{e}_{j_2}={\rm O}_p\left(\frac{1}{\sqrt n}\right),
\end{gather*}
where $\frac{1}{n}{\rm E}\left\{\text{tr}\mathbf{M}^{-2}(\theta_{i})\right\}\rightarrow y$.
It follows that
\begin{gather*}
\max_{1\leq j\leq q}\left|\frac{1}{n}\textbf{e}_j^{\top}(1+\delta_{i})\mathbf{Z}_{1}\mathbf{Z}_{2}^{\top}\mathbf{M}^{-2}(\theta_{i})\left(\frac{1}{n}\mathbf{Z}_{2}\mathbf{Z}_{2}^{\top}\right)^{-1}\frac{1}{n}\mathbf{Z}_{2} \mathbf{Z}_{1}^{\top}\textbf{e}_j-\frac{1}{n}{\rm E}\left[\text{tr}\left\{\mathbf{M}^{-2}(\theta_{i})\right\}\right]\right|={\rm O}_p\left(\frac{\sqrt q}{\sqrt n}\right)\\
\max_{1\leq j_1\neq j_2\leq q}\left|\frac{1}{n}\textbf{e}_{j_1}^{\top}(1+\delta_{i})\mathbf{Z}_{1}\mathbf{Z}_{2}^{\top}\mathbf{M}^{-2}(\theta_{i})\left(\frac{1}{n}\mathbf{Z}_{2}\mathbf{Z}_{2}^{\top}\right)^{-1}\frac{1}{n}\mathbf{Z}_{2} \mathbf{Z}_{1}^{\top}\textbf{e}_{j_2}\right|={\rm O}_p\left(\frac{q}{\sqrt n}\right).
\end{gather*}
Similarly, we can get the following for other terms:
\begin{gather*}
\max_{1\leq j\leq q}\left|\textbf{e}_j^{\top}(1+\delta_{i})\frac{1}{n} \mathbf{Z}_{1}\mathbf{Z}_{1}^{\top}\textbf{e}_j-1\right|={\rm O}_p\left(\frac{\sqrt q}{\sqrt n}\right),\
\max_{1\leq j_1\neq j_2\leq q}\left|\textbf{e}_{j_1}^{\top}(1+\delta_{i})\frac{1}{n} \mathbf{Z}_{1}\mathbf{Z}_{1}^{\top}\textbf{e}_{j_2}\right|={\rm O}_p\left(\frac{ q}{\sqrt n}\right),\\
\max_{1\leq j\leq q}\left|\textbf{e}_j^{\top}\frac{1}{\hat\lambda_{i} \theta_{i}}
\frac{1}{T}\mathbf{X}_{1}\mathbf{X}_{2}^{\top}\mathbf{M}^{-1}(\hat\lambda_{i})\mathbf{M}^{-1}(\theta_{i})\left(\frac{1}{n}\mathbf{Z}_{2}\mathbf{Z}_{2}^{\top}\right)^{-1}\frac{1}{T}\mathbf{X}_{2}
\mathbf{X}_{1}^{\top}\textbf{e}_j\right|={\rm O}_p\left(\frac{\sqrt{\sum_j\lambda_j^2}}{\lambda_i^2}\right),\\
\max_{1\leq j_1\neq j_2\leq q}\left|\textbf{e}_{j_1}^{\top}\frac{1}{\hat\lambda_{i} \theta_{i}}
\frac{1}{T}\mathbf{X}_{1}\mathbf{X}_{2}^{\top}\mathbf{M}^{-1}(\hat\lambda_{i})\mathbf{M}^{-1}(\theta_{i})\left(\frac{1}{n}\mathbf{Z}_{2}\mathbf{Z}_{2}^{\top}\right)^{-1}\frac{1}{T}\mathbf{X}_{2}
\mathbf{X}_{1}^{\top}\textbf{e}_{j_2}\right|={\rm O}_p\left(\frac{{\sum_j\lambda_j}}{\lambda_i^2}\right),\\
\max_{1\leq j\leq q}\left|\textbf{e}_j^{\top}\frac{(1+\delta_{i})}{\hat\lambda_{i}} \frac{1}{n} \mathbf{Z}_{1}
\mathbf{Z}_{2}^{\top}\mathbf{M}^{-1}(\hat\lambda_{i})\mathbf{M}^{-1}(\theta_{i})\left(\frac{1}{n}\mathbf{Z}_{2}\mathbf{Z}_{2}^{\top}\right)^{-1}\frac{1}{T}\mathbf{X}_{2}\mathbf{X}_{1}^{\top}  \textbf{e}_j\right|={\rm O}_p\left(\frac{\sqrt{\sum_j\lambda_j}}{\lambda_i}\right),\\
\max_{1\leq j_1\neq j_2\leq q}\left|\textbf{e}_{j_1}^{\top}\frac{(1+\delta_{i})}{\hat\lambda_{i}} \frac{1}{n} \mathbf{Z}_{1}
\mathbf{Z}_{2}^{\top}\mathbf{M}^{-1}(\hat\lambda_{i})\mathbf{M}^{-1}(\theta_{i})\left(\frac{1}{n}\mathbf{Z}_{2}\mathbf{Z}_{2}^{\top}\right)^{-1}\frac{1}{T}\mathbf{X}_{2}\mathbf{X}_{1}^{\top}\textbf{e}_{j_2}\right|
={\rm O}_p\left(\frac{\sqrt{q\sum_j\lambda_j}}{\lambda_i}\right),\\
\max_{1\leq j\leq q}\left| \textbf{e}_j^{\top}\frac{1}{n} \mathbf{Z}_{1} \mathbf{C}(\theta_{i})\mathbf{X}_{1}^{\top}\textbf{e}_j\right|={\rm O}_p\left(\frac{\sqrt{\sum_j\lambda_j}}{\sqrt {n}\lambda_i}\right),\
\max_{1\leq j_1\neq j_2\leq q}\left|\textbf{e}_{j_1}^{\top}\frac{1}{n} \mathbf{Z}_{1} \mathbf{C}(\theta_{i})\mathbf{X}_{1}^{\top}\textbf{e}_{j_2}\right|={\rm O}_p\left(\frac{\sqrt{q\sum_j\lambda_j}}{\sqrt {n}\lambda_i}\right),\\
\max_{1\leq j\leq q}\left|\textbf{e}_j^{\top}\frac{(1+\delta_{i})}{\hat\lambda_{i}} \frac{1}{T} \mathbf{X}_{1}
\mathbf{X}_{2}^{\top}\mathbf{M}^{-1}(\hat\lambda_{i})\mathbf{M}^{-1}(\theta_{i})\left(\frac{1}{n}\mathbf{Z}_{2}\mathbf{Z}_{2}^{\top}\right)^{-1}\frac{1}{n}\mathbf{Z}_{2}\mathbf{Z}_{1}^{\top}\textbf{e}_j\right|={\rm O}_p\left(\frac{\sqrt{\sum_j\lambda_j}}{\lambda_i}\right),\\
\max_{1\leq j_1\neq j_2\leq q}\left|\textbf{e}_{j_1}^{\top}\frac{(1+\delta_{i})}{\hat\lambda_{i}} \frac{1}{T} \mathbf{X}_{1}
\mathbf{X}_{2}^{\top}\mathbf{M}^{-1}(\hat\lambda_{i})\mathbf{M}^{-1}(\theta_{i})\left(\frac{1}{n}\mathbf{Z}_{2}\mathbf{Z}_{2}^{\top}\right)^{-1}\frac{1}{n}\mathbf{Z}_{2}\mathbf{Z}_{1}^{\top}\textbf{e}_{j_2}\right|={\rm O}_p\left(\frac{\sqrt{q\sum_j\lambda_j}}{\lambda_i}\right),\\
\max_{1\leq j\leq q}\left|\textbf{e}_j^{\top}\frac{1}{T}\mathbf{X}_{1} \mathbf{D}(\theta_{i})\mathbf{Z}_{1}^{\top}\textbf{e}_j\right|={\rm O}_p\left(\frac{\sqrt{\sum_j\lambda_j}}{\sqrt {n}\lambda_i}\right),\
\max_{1\leq j_1\neq j_2\leq q}\left|\textbf{e}_{j_1}^{\top}\frac{1}{T}\mathbf{X}_{1} \mathbf{D}(\theta_{i})\mathbf{Z}_{1}^{\top}\textbf{e}_{j_2}\right|={\rm O}_p\left(\frac{\sqrt{q\sum_j\lambda_j}}{\sqrt {n}\lambda_i}\right).
\end{gather*}

Thus, all these inequalities  lead to
\begin{gather*}
\max_{1\leq j\leq q}\left|\textbf{e}_j^{\top}\Theta_{2n}\textbf{e}_j-(y-1)\right|={\rm O}_p\left(\frac{\sqrt q \delta_i}{\lambda_i}+\frac{\sqrt q}{\sqrt n}+\frac{\sqrt{\sum_j\lambda_j^2}}{\lambda_i^2}+\frac{\sqrt{\sum_j\lambda_j}}{\lambda_i}\right),\\
\max_{1\leq j_1\neq j_2\leq q}\left|\textbf{e}_{j_1}^{\top}\Theta_{2n}\textbf{e}_{j_2}\right|={\rm O}_p\left(\frac{ q\delta_i}{\lambda_i}+\frac{q}{\sqrt n}+\frac{{\sum_j\lambda_j}}{\lambda_i^2}+\frac{\sqrt{q\sum_j\lambda_j}}{\lambda_i}\right).
\end{gather*}
The proof is completed. \quad$\square$

\section*{Acknowledgments}
The authors gratefully acknowledge a grant from the University Grants Council of Hong Kong and a NSFC grant (NSFC11671042). Drs Xie and Zeng are co-first authors. 

\end{document}